\documentclass[12pt]{amsart}
\usepackage{amsmath}\usepackage{MnSymbol}
\usepackage{skak}
\usepackage{amsfonts}\usepackage{verbatim}\usepackage{bbold}
\usepackage{graphicx}\usepackage{fontawesome}
\pdfoutput=1
\usepackage{marvosym}\usepackage{ifsym}
\usepackage{dictsym}\usepackage{wasysym}
\usepackage{pigpen}\pdfgentounicode=0
\usepackage{amstext}\usepackage{bclogo}
\usepackage{amsbsy}\usepackage{tikzsymbols}
\usepackage{amsopn}\usepackage{hyperref}%
\usepackage{thmtools}
\usepackage{amsthm}\usepackage{color}

\def\proclaim#1{\vskip0.5em\noindent{\bf #1}\it }
\def\endproclaim{\vskip0.5em\par\noindent\rm}

\def\proclaim#1{\vskip0.5em\noindent{\bf #1}\it}
\def\endproclaim{\vskip0.5em\par\noindent\rm}

\def\undersetbrace#1\to#2{\underbrace{#2}_{#1}}

\def\demo#1{\vskip0.5em\noindent{\bf #1\ }}

\def\text#1{\mbox{#1}}
\def\flushpar{\par\noindent}

\def\tag#1{\eqno{(#1)}}

\newcommand{\mapright}[1]{%
    \smash{\mathop{%
        \hbox to 1cm{\rightarrowfill}
        }
    \limits^{#1}
    }
}
\newcommand{\mapleft}[1]{%
    \smash{\mathop{%
        \hbox to 1cm{\rightarrowfill}
        }
    \limits_{#1}
    }
}

\def\e{\epsilon}
\def\a{\alpha}

\def\G{\Gamma}
\def\g{\gamma}
\def\d{\delta}
\def\D{\Delta}
\def\s{\sigma}

\def\th{\theta}
\def\l{\lambda}

\def\x{\times}

\def\f{\flushpar}
\def\u{\underline}
\def\v{\varphi}

\def\om{\omega}
\def\Om{\Omega}
\def\B{\mathcal B}
\def\T{\widehat T}
\def\({\bigg(}
\def\){\bigg)}
\def\<{\langle}
\def\>{\rangle}

\def\bul{\smallskip\f$\bullet\ \ \ $}\def\lfl{\lfloor}\def\rfl{\rfloor}

\def\with respect to{\text{with respect to}}

\def\<{\bold\langle}
\def\>{\bold\rangle}
\def\pprime{\prime\prime}

\def\bul{\smallskip\f$\bullet\ \ \ $}\def\sms{\smallskip\f}\def\sbul{\f$\bullet\
\ \ $}\def\Par{\smallskip\f\P}

\def\ttau{\widehat{\tau}}
\begin{document}

    \title{Functional limits for ``tied down''  occupation time processes of infinite ergodic transformations}
\author{ Jon. Aaronson and Toru Sera}
\address[Aaronson]{School of Math. Sciences, Tel Aviv University,
69978 Tel Aviv, Israel.}
\email{aaro@tau.ac.il}
\address[Sera]{Department of Mathematics, Graduate School of Science, Osaka University,
Toyonaka, Osaka 560-0043, Japan}
\email{sera@math.sci.osaka-u.ac.jp}
\begin{abstract}  We prove functional, distributional limit theorems  for the occupation times  of
pointwise dual ergodic transformations at  ``tied-down'' times immediately after ``excursions''.
The limiting processes are  tied down Mittag-Leffler processes
 and
the transformations involved exhibit  functional versions of the tied-down renewal
properties in \cite{A-Sera}.
\end{abstract}
\subjclass[2010]{37A40, 60F05}
\keywords{infinite ergodic theory,  pointwise dual ergodic, Darling-Kac theorem, Mittag Leffler processes,
 stable subordinator, rational weak mixing, local  time  process, conditional limit theorem, tied-down renewal theory.}
\thanks{\copyright 2019-21. Aaronson{\tiny$^\prime$}s research was partially supported by ISF grant No. 1289/17. Sera{\tiny$^\prime$}s research was partially supported by JSPS KAKENHI grants No.\ JP19J11798 and No.\ JP21J00015.}
\maketitle\markboth{Functional limits for ``tied down''  occupation time processes}{Jon. Aaronson and Toru Sera}

\section*{\S0 Introduction}
\

\subsection*{Overview}
\

\ \  For a  {\tt\small pointwise dual ergodic}  measure preserving transformation (see p.\pageref{PDE})  
with  $\g$-regularly varying return sequence 
 ($\g\in (0,1)$) the $\g$-{\tt\small Mittag-Leffler process} ($\frak m_\g$ on p.\pageref{ML})  
 appears as a scaling limit of   the normalized Birkhoff sums of non-negative, integrable functions (\ref{dsrailways} on p.\pageref{dsrailways}).
\
 
 The {\tt\small tied down $\g$-Mittag-Leffler process} $\frak w_\g$
 is a functional of  $\frak m_\g$,  defined by \ref{Football} on p. \pageref{Football}  \ so that, in particular, $1$ is a point of increase (i.e.  a.s. $\frak w_\g(t)<\frak w_\g(1)\ \forall\ t<1$). 
 \
 
 In this paper we consider similar functionals of
  the Birkhoff sum processes:  conditioning  on arrival under the increment semiflow (p. \pageref{increment}) in a fixed, finite measure set at the final observation time. In particular, we  show that their weighted
 scaling limits, conditioned thus,  are tied-down Mittag-Leffler processes  (Theorem U on p.\pageref{thmU}).
\

This kind of conditioning has been  studied in the context of  Markov  processes  (e.g. \cite{Dw-Ka,Wendel,Liggett1968,BPY,FPY,Ber,ChaUri}).
The tied-down Mittag-Leffler processes appear as {\tt\small local times} of {\tt\small Markov bridges} (see \S6).

\

 \subsection*{Non-decreasing stochastic processes}

\

Let 
$ D_{\uparrow,R}$ denote the collection of distribution functions of finite Borel measures on $[0,R]$ ($R\in\Bbb R_+$) and let $ D_{\uparrow,\infty}$ denote the distribution functions of Radon measures on $[0,\infty)$:
\begin{align*}
 D_{\uparrow,R}:=\{\xi :[0,R]\to [0,\infty):\ \xi\ \text{ is non-decreasing},\ \xi (x+)=\xi (x)\ \forall\ x\in [0,R]\}
\end{align*}
where   $\xi (x+):=\lim_{y\to x,\ y>x}\xi (y)$ ($0\le x<\infty$)  $\&\ \xi (x-):=\lim_{y\to x,\ y<x}\xi (y)$   ($0< x\le\infty$);
and let
$$D_{\uparrow,\infty}:=\{\xi :[0,\infty)\to [0,\infty):\ \xi |_{[0,R]}\in D_{\uparrow,R}\ \forall\ R>0\}.$$
\

Equip $D_{\uparrow,R}$ ($0<R<\infty$) with the $\s$-algebra
$$\B_R:=\s(\{e_t:\ 0\le t\le R\})\ \text{where for}\ \xi\in D_{\uparrow,R},\ e_t(\xi):=\xi(t)$$
and equip $D_{\uparrow,\infty}$ with $\B_\infty:=\s(\{e_t:\ t\in [0,\infty)\})$.
\

A {\it non-decreasing stochastic process} on $[0,R]$ is a random variable $Z$ on the measurable space
$(D_{\uparrow,R},\B_R)$.  It    is  {\it a.s. continuous} if a.s., $Z$ is continuous; and {\it continuous in probability} if
$\forall\ t\in [0,R]$, $Z(s)\rightarrow Z(t)$ as $s\to t$ in probability;
equivalently $Z(t-)=Z(t)$ a.s. $\forall\ t\in [0,R]$.

\

\subsection*{Convergence}
\

The {\it Skorokhod $J_1$ topology }  on $D_{\uparrow,R}$ ($R\le\infty$) 
is generated by the Polish metric $J_1^{(R)}$ (\cite{Skorokhod, E-K, Billingsley}) defined   for $R<\infty$ by
\begin{align*}
 J_1^{(R)}(\xi,\eta):=\inf\,\{\|\xi&-\eta\circ\ell\|_{L^\infty([0,R])}+\|\log(\ell')\|_{L^\infty([0,R])}:\ \ell:[0,R]\to [0,R]\\ & \text{an increasing bi-Lipschitz homeomorphism}\ \}
\end{align*}
and for  $R=\infty$ by 
\begin{align*}
 J_1^{(\infty)}(&\xi,\eta):=\inf\{\|\log(\ell')\|_{L^\infty([0,\infty))}+
 \int_0^\infty e^{-u}(\|\xi^{(u)}-\eta^{(u)}\circ\ell\|_{L^\infty([0,u])}\wedge 1 )du:\\ & 
 \ell:[0,\infty)\to [0,\infty)\ \text{an increasing bi-Lipschitz homeomorphism}\ \}
\end{align*}
where $x^{(u)}(t):=x(t\wedge u)$. 
\

It follows that for $\xi_n,\ \xi\in D_{\uparrow,\infty}$,
\

$$\xi_n\xrightarrow[n\to\infty]{J_1^{(\infty)}}\xi\ \Leftrightarrow\ \xi_n^{(u)}\xrightarrow[n\to\infty]{J_1^{(u)}}\xi^{(u)}\ \forall\ \text{continuity points $u$ of $\xi$}.$$

The {\it distribution function} ({\tt DF}) topology on  $D_{\uparrow,R}$ 
is that inherited from the weak $*$ topology on $C([0,R])^*$ for $R<\infty$  
and from the weak $*$ topology on $C_C([0,\infty))^*$ for $R=\infty$ where $C_C([0,\infty))$ denotes the space of continuous functions with compact support.

The spaces  $(D_{\uparrow,R},{\tt DF})$ are Polish spaces and locally compact for   $R<\infty$.
\

Throughout, given a Polish space $Y$,
\bul $C_B(Y)$ is  the space of bounded, continuous, $\Bbb R$-valued functions
equipped with the norm $\|g\|_{C_B}:=\sup_{y\in Y}|g(y)|$,
\

and
\bul $\B(Y):=\s(\{\text{open sets}\})=\s(C_B(Y))$ is  the Borel $\s$-algebra.

\subsection*{Comparison: {\tt DF} vs $J_1$ topologies}\label{Comparison}
\

 {\tt DF} convergence is pointwise convergence at all continuity points of the limit and
 follows from $J_1$ convergence.
\

{\tt DF} convergence to  a  continuous  limit entails (by monotonicity) uniform convergence
on bounded subintervals, whence $J_1$ convergence.
\

There are sequences converging {\tt DF} (to discontinuous limits) which do not converge $J_1$.

\

As shown in \cite{Skorokhod}, 
$$\B(D_{\uparrow,R},J_1)=\B(D_{\uparrow,R},{\tt DF})=\B_R.$$
\subsection*{Distributional convergence}\label{dist}
\

Given a Polish space $Y$, let
$$\text{\tt RV}\,(Y):=\{\text{\tt random variables on}\ Y\}$$
equipped with the topology of {\it weak convergence}  so that
$$Z_n\xrightarrow[n\to\infty]{\text{\tt\tiny RV}\,(Y)}Z\ \text{if}\ \Bbb E(g(Z_n))\xrightarrow[n\to\infty]{}\Bbb E(g(Z))
\ \forall\ g\in C_B(Y).$$
For $Z\in\text{\tt RV}\,(Y)$ the {\it distribution} of $Z$ is that probability measure
$\text{\tt dist}\,Z\in\mathcal{P}(Y)$ so that
$$\Bbb E(g(Z))=\int_Ygd(\text{\tt dist}\,Z)\ \text{for}\  g\in C_B(Y).$$

Note that for $R< \infty$, a collection $\mathcal Z\subset \text{\tt RV}\,(D_{\uparrow,R},{\tt DF})$ is {\it tight} (i.e. precompact) in $\text{\tt RV}\,(D_{\uparrow,R},{\tt DF})$ iff 
$\{Z(R):\ Z\in\mathcal{Z}\}$ is precompact in $\text{\tt RV}\,(\Bbb R_+)$ ($R<\infty$).
\

A collection $\mathcal Z\subset \text{\tt RV}\,(D_{\uparrow,\infty},{\tt DF})$ is tight when  $\{Z(r):\ Z\in\mathcal{Z}\}$ is precompact in $\text{\tt RV}\,(\Bbb R_+)\ \forall\ r>0$.
\

Convergence in the space $\text{\tt RV}\,(D_{\uparrow,R},{\tt DF})$ is characterized by convergence of finite dimensional distributions ({\tt FDD}s) (see \cite{Kal02}): 
\begin{align*}\tag*{{\tt FDD}}\label{FDD}&Z_n\xrightarrow[n\to\infty]{\text{\tt\tiny RV}\,(D_{\uparrow,R},{\tt DF})}\ Z\ \iff\ \ Z_n\xrightarrow[n\to\infty]{\text{\tt\tiny FDD}}\ Z
\ \ \ \text{\tt\large i.e.}:\\ &
(Z_n(t_1),Z_n(t_2),\dots,Z_n(t_k),Z_n(R))\xrightarrow[n\to\infty]{\text{\tt\tiny RV}\,(\Bbb R_+^{k+1})}\ (Z(t_1),Z(t_2),\dots,Z(t_k),Z(R))\\ & 
\ \ \ \ \ \ \ \ \  \forall\ k\ge 1,\ (t_1,t_2,\dots,t_k)\in (0,R)^k\ \text{with}\ \Bbb P(Z(t_i)=Z(t_i-))=1;
\end{align*}
and $Z_n\xrightarrow[n\to\infty]{\text{\tt\tiny RV}\,(D_{\uparrow,\infty},{\tt DF})}\ Z$\ \ iff
  $Z_n|_{[0.R]}\xrightarrow[n\to\infty]{\text{\tt\tiny RV}\,(D_{\uparrow,R},{\tt DF})}\ Z|_{[0.R]}$
  $\forall\ R\in\Bbb R_+$ with $\Bbb P(Z(R)=Z(R-))=1$.

\subsection*{Distributional convergence to an a.s. continuous limit}\label{aslt}
\

It follows from the comparison remarks above (p. \pageref{Comparison}) that 
\bul  convergence in $\text{\tt RV}\,(D_{\uparrow,R}, J_1)$ implies {\tt FDD} convergence; 
\

and,  
\bul {\tt FDD} convergence to an a.s. continuous  limit $Z$ with $Z(0)=0$ a.s. is sufficient for convergence  in $\text{\tt RV}\,(D_{\uparrow,R}, J_1)$ (\cite[Thm. VI.3.37]{Jacod-Shiryaev}). \footnote{\tiny However
convergence in $\text{\tt RV}\,(D_{\uparrow,R}, J_1)$ is not in general characterized by convergence of {\tt FDD}s even when the limit is continuous in probability as shown by the example in \cite[Problem 12.5]{Billingsley}.}

\

The  test functions used here to prove {\tt FDD} convergence to an a.s. continuous  limit include:
\subsection*{Product functions}
 A {\it product function} is a function  $\mathbb{g}_{\u t}:D_{\uparrow,1}\to\Bbb R$ of form
\begin{align*}\tag*{\dsliterary}\label{dsliterary}
 \mathbb{g}_{\u t}(\xi)=\big[\prod_{\nu=1}^Ng_\nu(\xi(t_\nu))\big]h(\xi(1)-\xi(t_N)) 
\end{align*}

where $\u t=(t_1,t_2,\dots,t_N)$ with $0\le t_1<t_2<\dots<t_N<1$ and 
 and $g_1,\dots,g_N,h\in C([0,\infty])_+$. 
 \
 
 Denote $\Pi:=\{{\tt\small product\  functions}\}$.
 \
 
 It is standard to show that 
 $\Pi\subset C_{B,c}(D_{\uparrow,1})$\label{cbc}, the collection of bounded, measurable functions on $D_{\uparrow,1}$  continuous at each  $\xi\in D_{\uparrow,1}\cap C([0,1])$ (i.e. at each $\xi:[0,1]\to\Bbb R$ which is continuous and non-decreasing).
 \

Moreover, for $Z_n,\ Z\in\text{\tt RV}(D_{\uparrow,1})$,
 \begin{align*}\tag*{{\scriptsize\faAnchor}}\label{faAnchor}\Bbb E(&\mathbb{g}(Z_n))\xrightarrow[n\to\infty]{×}\Bbb E(\mathbb{g}(Z))\ \forall\ \mathbb{g}\in\Pi\\ &
 \implies\ Z_n\xrightarrow[n\to\infty]{\text{\tt\tiny FDD}}\ Z.
 \end{align*}

\subsection*{Self similarity, increment stationarity and their dynamics}\label{increment}
\ 

The space $D_{\uparrow,\infty}$ is naturally equipped with
\bul the {\it increment semiflow} $T_s:D_{\uparrow,\infty}\to D_{\uparrow,\infty}$ ($s\ge 0$) defined by
$T_sZ(t):=Z(s+t)-Z(s)$;

\

and
\

\bul for $\g>0$, the $\g$-{\it scalings} $\D_{a,\g}:D_{\uparrow,\infty}\to D_{\uparrow,\infty}$ ($a\ge 0$) defined by
$\D_{a,\g}Z:=\tfrac1{a^\g}Z\circ D_a$ with
$D_a(t):=at$.

\

For $\g> 0$, we call the  process $Z\in\text{\tt RV}\,(D_{\uparrow,\infty})$ $\g$-{\it self similar}  if 
$\D_{a,\g}Z\overset{d}=Z\ \forall\ a>0.$ 
\

Equivalently, $(D_{\uparrow,\infty},\text{\tt\small dist}\,Z,\D_{a,\g})$ is a probability preserving
transformation  $ \forall\ a>0$. 

\

Also, we'll say that the process $Z\in\text{\tt RV}\,(D_{\uparrow,\infty})$ {\it has stationary increments}  if
$T_sZ\overset{d}=Z\ \forall\ s>0$. 
\

Equivalently,  $(D_{\uparrow,\infty},\text{\tt\small dist}\,Z,T_s)$ is a probability preserving
transformation  $ \forall\ s>0$.

\

\subsection*{Inverse process}\label{inverse}
\

The {\it inverse} of $\xi\in D_{\uparrow,\infty}$ with $\xi(\infty-)=\infty$
is $\mathcal{I}(\xi):=\xi^{-1}\in D_{\uparrow,\infty}$ defined by
$$\xi^{-1}(t):=\inf\,\{s>0:\ \xi(s)> t\}.$$ 
\

If $Z\in\text{\tt RV}\,(D_{\uparrow,\infty})$ is $\g$-self similar with $\g>0$ and  $Z(\infty-)=\infty$ a.s., then $Z^{-1}$ is $\tfrac1\g$-self similar.
\subsection*{Waiting times}\label{waiting}
\

Define the {\it waiting time functionals}
$\mathcal{G},\ \mathcal{D}$ on $D_{\uparrow,\infty}$ as follows:
$$\mathcal G(\xi)(t):=\inf\,\{s\le t:\ \xi(s)=\xi(t)\}$$
so that $\mathcal G:D_{\uparrow,\infty}\to D_{\uparrow,\infty}$; and
$$\mathcal D(\xi)(t):=\sup\,\{s\ge t:\ \xi(s)=\xi(t)\}\le\infty.$$
If $\xi(\infty-)=\infty$, then $\mathcal D(\xi)\in D_{\uparrow,\infty}$.

\

\subsection*{Subordinators and Mittag-Leffler processes}
\

For $\g\in (0,1)$, 
\bul the $\g$-{\it stable subordinator} is $\eta_\g\in\text{\tt RV}\,( D_{\uparrow,\infty})$, $\tfrac1\g$-self similar with positive, stationary, independent $\g$-stable increments so that 
$$\Bbb E(\eta_\g^{-1}(1))=\Bbb E(\tfrac1{\eta_\g(1)^\g})=1;$$
\bul\label{ML} the $\g$-{\it Mittag-Leffler ({\tt $\g$-ML}) process} is  $\frak m_\g=\eta_\g^{-1}$ and is $\g$-self similar;
 
\bul  the {\it tied down $\g$-{\tt ML}  process} is $\frak w_\g\in\text{\tt RV}\,( D_{\uparrow,1})$ defined 
by
\begin{align*}\tag*{\Football}\label{Football}\frak w_\g(t):=\frac{\frak m_\g(g_{\g,1}t)}{g_{\g,1}^\g}=\D_{g_{\g,1},\g}\frak m_\g(t)
\end{align*}
where $g_{\g,1}:=\mathcal G(\frak m_\g)(1)$.

\

 These processes $\frak{w}_\g$ for $0<\g\le\frac12$ correspond to the local time at zero of the {\tt\small symmetric $\frac1{1-\g}$-stable bridge}; and for $0<\g<1$, to the local time at zero of the {\tt\small  Bessel bridge of dimension} $2-2\gamma$ (see \S6). 
 \subsection*{Almost sure continuity of the {\tt ML} processes}
 \
 
The stable subordinators $\eta_\g$ are strictly increasing a.s.;
 thus, a.s., $\frak m_\g=\eta_\g^{-1}$ is continuous (see e.g. the remark on p.8 of
 \cite{RevYor}); as is $\frak w_\g$  by \ref{Football}.
 
\subsection*{Pointwise dual ergodic transformations}
\

Let $(X,m,T)$   denote   a measure preserving transformation  $T$ of the non-atomic, Polish measure space $(X,m)$ (where 
$m$ is a $\s$-finite, non-atomic measure
defined  on the Borel subsets $\B(X)$ of $X$). 
\

The associated {\it transfer operator} $\T:L^1(m)\to L^1(m)$ is the predual of $f\mapsto f\circ T\ \ (f\in L^\infty(m))$, that is 
$$\int_X\T f gdm=\int_Xf g\circ Tdm\ \ \text{for}\ f\in L^1(m)\ \&\ g\in L^\infty(m).$$
\

The measure preserving
transformation  $(X,m,T)$ is called  {\it pointwise dual ergodic} if there is a sequence $a(n)=a_n(T)$ (the {\it return sequence} of $(X,m,T)$) so that
\begin{align*}\tag*{{\tt PDE}}\label{PDE}\frac1{a(n)}\sum_{k=0}^{n-1}\T^kf\xrightarrow[n\to\infty]{}\ \int_Xfdm\ \text{a.e.}\ \forall\ f\in L^1(m). 
\end{align*}

Pointwise dual ergodicity entails 
\sbul {\tt conservativity} (aka  recurrence-- no non-trivial wandering sets); 

\sbul{\tt ergodicity} (no non-trivial invariant sets) and 

\sbul {\tt rational ergodicity} as defined in \cite[\S 3.8]{IET} (see also \ref{faTree} on p. \pageref{faTree}).

\subsection*{Stationary processes and skyscrapers}
\

A (discrete) {\it stationary process}   is a quadruple $(\Om,\mu,\tau,\phi)$ where $(\Om,\mu,\tau)$ is a probability preserving
transformation  and $\phi:\Om\to\Bbb X$ (some metric space) is measurable. The stationary process $(\Om,\mu,\tau,\phi)$ is called {\it ergodic} if 
$(\Om,\mu,\tau)$ is an ergodic  probability preserving
transformation 

For example if $(X,m,T)$ is a conservative measure preserving
transformation  and 

$\Om\in\mathcal F_+:=\{A\in \B(X):\ \ 0<m(A)<\infty\}$\label{mathcalF},\ the {\it return time function} to $\Om$ is $\v=\v_\Om:\Om\to\mathbb N$ defined by
$\v(\om):=\min\{n\ge 1:\ T^n\om\in\Om\}<\infty$ a.s. by conservativity. 

The {\it induced transformation} on  $\Om$ is $T_\Om:\Om\to\Om$ defined by $T_\Om(\om):=T^{\v(\om)}(\om)$. As shown in
\cite{Kakutani}, $(\Om,m_\Om,T_\Om)$ is a probability preserving
transformation  which, in case $X\overset{m}=\bigcup_{n\ge 0}T^{-n}\Om$, is ergodic together with $(X,m,T)$. The $\Bbb N$-valued stationary process $(\Om,m_\Om,T_\Om,\v_\Om)$
is called the {\it return time process} to $\Om$.
\

As in \cite{Kakutani}, the {\it skyscraper}  over the $\Bbb N$-valued stationary process $(\Om,\mu,\tau,\phi)$ is the conservative measure preserving
transformation 
 $(X,m,T)=(\Om,\mu,\tau)^\phi$ defined by
 \begin{align*}
 X:=\{(\om,n)\in\Om\x\Bbb N:\ &1\le n\le \phi(\om)\},\ m:=\mu\x\#|_X\ \&\\ & T(\om,n):=\begin{cases}& (\om,n+1)\ \ n<\phi(\om)\\ &
                                                                                   (\tau(\om),1)\ \ n=\phi(\om).
                                                                                  \end{cases}
\end{align*}

\

If $(X,m,T)$ is a conservative, ergodic,  measure preserving 
transformation , $\Om\in\B(X),
\ 0<m(\Om)<\infty$, then $(X,m,T)$ is a  factor of $(\Om,m_\Om,T_\Om)^{\v_\Om}$
    (isomorphic if $(X,m,T)$ is  invertible) and
    the return time process of $(\Om,m_\Om,T_\Om)^{\v_\Om}$ to  $\Om\x\{1\}$  is isomorphic with $(\Om,\mu,\tau,\phi).$ 
\subsection*{Functional distributional limits}\ \ Let $(X,m,T)$ be a pointwise dual ergodic measure preserving
transformation  with $\g$-regularly varying return sequence
$a(n)=a_n(T)$ ($0<\g<1$).
\

By the functional Darling-Kac Theorem (\cite{Bin71,Ow-Sam}),
\begin{align*}\tag*{\dsrailways}\label{dsrailways}
\Psi_{n,f}\xrightarrow[n\to\infty]{\frak d}\ m(f)\frak m_\g\ \text{in}\ (D_{\uparrow,\infty},J_1)\ \forall\ f\in L^1_+:=\{g\in L^1(m):\ g\ge 0,\ \int_Xgdm>0\} 
\end{align*}

 where $\Psi_{n,f}(t):=\tfrac{S_{\lfl nt\rfl}(f)}{a(n)}$ with $S_n(f):=\sum_{k=0}^{n-1}f\circ T^k$ is the {\it normalized Birkhoff-sum step function} and where ``$\xrightarrow{\frak d}$ in $\Xi$'' means distributional convergence in the metric space $\Xi$ with respect to all
$m$-absolutely continuous probabilities as obtained by Eagleson's theorem (\cite{Eagleson, TZ}).
\

By \cite{TK-cgce, TK-weak}, if $(\Om,\mu,\tau,\phi)$ is an exponentially continued fraction mixing, $\Bbb R_+$-valued stationary process, that is
\par for some $K'>0,\ \th\in (0,1)$,
$$|\mu(a\cap \tau^{-(n+k)}B)-\mu(a)\mu(B)|\le K'\th^n\mu(a)\mu(B)\ \ \forall\
n,k\ge 1,\ a\in\s(\{\phi\circ\tau^j\}_0^{k-1}),\ B\in\B(\Om);$$

and $\frac1{\mu([\phi>t])}$ is $\g$-regularly varying with $0<\g<1$, then 
with
$$a(t):=\frac1{\G(1-\g)\G(1+\g)
\mu([\phi>t])},$$ 
\begin{align*} \tag*{{\scriptsize\faMedkit}}\label{faMedkit}
 \Phi_n\xrightarrow[n\to\infty]{\frak d}\eta_\g \ \text{ in } \ (D_{\uparrow,\infty},J_1)
\end{align*}
 where
$\Phi_n(t):=\tfrac{\phi_{\lfl nt\rfl}}{a^{-1}(n)}$ with $\phi_k:=\sum_{j=0}^{k-1}\phi\circ\tau^j$.
\

In the sequel, we'll prove tied down versions of (\dsrailways) for conservative, ergodic,  measure preserving 
transformations having (among other properties) a return time process satisfying ({\scriptsize\faMedkit}).

\section*{\S1 Results}
\subsection*{Tied down Mittag-Leffler processes}
\

The following propositions give characterizations of  $\frak w_\g$ defined by \ref{Football} (on p. \pageref{Football}):
\proclaim{Proposition C\ \ {\rm (characterization of  $\frak w_\g$)}}\label{propC}
\

\ For $\g\in (0,1)$, 

\begin{align*}\tag*{\symbishop}\label{bishop}
\Bbb E(h(\frak w_\g))=\Bbb E(\int_0^1&h(\D_{t,\g}\frak m_\g)d\frak m_\g(t))\\ & \forall\ h\in C_{B}( D_{\uparrow,1}):=C_{B}( D_{\uparrow,1},J_1).\end{align*}
\endproclaim
It follows by choosing $h(\xi)=H(\xi(1))$ ($H\in C_B(\Bbb R)$) in \ref{bishop}  that
$$\Bbb E(H(\frak w_\g(1)))=\Bbb E(\frak m_\g H(\frak m_\g(1)))\ \ \text{for}\ H\in C_B(\Bbb R).$$
\proclaim{Proposition D  {\rm ({\tt FDD}s of  $\frak w_\g$)} }\label{propD}
\

    \ \ For $\mathbb{g}_{\u t}:D_{\uparrow,1}\to\Bbb R_+$ a  product function \footnote{\small\rm see \ref{dsliterary} on p. \pageref{dsliterary}} 
 of form
\

$\mathbb{g}_{\u t}(\xi)=\frak g(\xi)h(\xi(1)-\xi(t_N))$
\ where $\frak g(\xi):=\prod_{\nu=1}^Ng_\nu(\xi(t_\nu))$
\

with
$0\le t_1<t_2<\dots<t_N<1$ and $g_1,\dots,g_N,h\in C([0,\infty])_+$,
{\small\begin{align*}\tag*{\symqueen}\label{queen} \Bbb E(&\mathbb{g}_{\u t}(\mathfrak w_\g))=\frak e_\g(\mathbb{g}_{\u t})\\ &:= \Bbb E\biggl(\frak g(\frak m_\g)1_{[\mathcal D(\frak m_\g)(t_N)\le 1]}(\tfrac1{1-\mathcal D(\frak m_\g)(t_N)})^{1-\g}h\biggl((1-\mathcal D(\frak m_\g)(t_N))^\g W_\g \biggr)\biggr)
\end{align*}}where $W_\g\perp \frak m_\g$ and $W_\g\overset{d}=\frak w_\g({\color{red} 1})$.
\endproclaim
\subsection*{Examples: Renewal shifts}
\

Write $\mathcal P(\Bbb N):=\{\text{\tt\small probability measures on}\ \Bbb N\}$.
\

For $f\in\mathcal{P}(\Bbb N)$, the associated Bernoulli  map is 
$$(\Om,p,S):=(\Bbb N^\Bbb N,f^\Bbb N,\text{\tt shift}).$$ 
  The {\it renewal shift}
corresponding to $f$ is the Kakutani skyscraper
$$(X,m,T):=(\Om,p,S)^\phi$$ where $\phi:\Om\to\Bbb N$ is defined by $\phi(x):=x_1$.
\

The {\it base} of the skyscraper  $(\Om,p,S)^\phi$ is $\widetilde{\Om}:=\Om\x\{1\}\in\B(X)$. It   is aka the {\it recurrent event} of the renewal shift and the corresponding distribution $f$ is known as its {\it lifetime distribution}.

\

The {\it renewal sequence} $(u(n):\ n\ge 0)$ corresponding  to $f\in\mathcal{P}(\Bbb N)$ is 
                                                                                            
$u(n)=u_f(n):=m(\widetilde{\Om}\cap T^{-n}\widetilde{\Om})$ and satisfies 
$u(n)=\sum_{k=1}^nf_ku(n-k)$ where (here and throughout) $f_k:=f(\{k\})$.
\

The renewal shift $(X,m,T)$ is a conservative, ergodic Markov shift, whence pointwise dual ergodic with return sequence  $a_n(T)=a(n)=a_f(n):=\sum_{k=1}^nu(k)$.
\

\

Now let $\g\in (0,1)$ and let $f\in\mathcal P(\Bbb N)$, then by Karamata's theorems
$$c(n):=p([\phi\ge n])=\sum_{k\ge n}f_k$$ is $(-\g)$-regularly varying iff $a(n)$ is $\g$-regularly varying and in this case,
$a(n)\propto\tfrac1{c(n)}$.\footnote{\label{propto}
Here and throughout, for $a_n,b_n>0,\ a_n\propto b_n$ means $\tfrac{a_n}{b_n}\xrightarrow[n\to\infty]{×}c\in\Bbb R_+$.}
\

\proclaim{ The Strong Renewal Theorem}
\ 

Let $f\in\mathcal P(\Bbb N)$ and let $a(n)$ be $\g$-regularly varying with $\g\in (0,1)$.

\sms{\rm (i)    (\cite{G-L,C-D})} If $f$ is non-arithmetic and $a(n)\gg\sqrt n$
\ \ \footnote{\label{ll}
Here and throughout, for $a_n,b_n>0,\ a_n\ll b_n$ means $\varlimsup_{n\to\infty}\tfrac{a_n}{b_n}<\infty$ and $a_n\sim b_n$ means $\tfrac{a_n}{b_n}\xrightarrow[n\to\infty]{×}1$.}, then
\begin{align*}\tag{{\tt SRT}}\label{SRT}
 u(n)\ \ \sim\ \ \frac{\g a(n)}n.
\end{align*}
\sms{\rm (ii)   (\cite{C-D})}\ 
For each $\g$-regularly varying  $a(n)\ngg\sqrt n$, $\exists\ f\in\mathcal P(\Bbb N)$ 
 aperiodic with $a_f(n)\sim a(n)$ and for which   {\rm({\tt SRT})}  fails.

\sms{\rm (iii)  (\cite[Thm. B]{Doney97} -- see also \cite{Goue11})}
\

If $f$ is  aperiodic and $f_n\ll\tfrac{c(n)}n$, then {\rm(\ref{SRT})} holds.
\endproclaim

For $(X,m,T)=(\Om,p,S)^\phi$  a renewal shift,  
write $$s_n:=S_n(1_{\widetilde{\Om}})\ \&\ \psi_n(t):=\tfrac{s_{\lfl nt\rfl}}{a(n)}.$$
\proclaim{Proposition E}\label{propE} 
\

\ \ Let $(X,m,T)=(\Om,p,S)^\phi$ be a renewal shift with aperiodic lifetime distribution
 $f\in\mathcal P(\Bbb N)$ so that
$a(n)$ is $\g$-regularly varying with $\g\in (0,1)$ and 
satisfying {\rm(\ref{SRT})},
 then 
 for $g\in C_B(D_{\uparrow,1})$,
\begin{align*}\tag*{\faBug}\label{faBug}
 \tfrac1{u(N)}\int_{\widetilde{\Om}\cap T^{-N}\widetilde{\Om}} g(\psi_{N})dm\xrightarrow[N\to\infty]{}\Bbb E(g(\frak w_\g)).
\end{align*}
\endproclaim
Proposition E is a strengthening of the Strong Renewal Theorem   as in \cite{G-L, Doney97,Goue11, C-D}
and is a special case of Theorem U\ref{Industry} (below on p.\pageref{Industry}).

\

For the general, renewal shift with regularly varying $a(n)$, we have weighted Cesaro versions of \ref{faBug} as in Theorem U, \ref{umbrella} $\&$ \ref{knight}.

\subsection*{Gibbs-Markov sets and cylinders}\label{gmcyl}
\

A {\it Markov map} is a quapruple $(Y,\nu,\tau,\a)$ where $(Y,\nu,\tau)$ is a nonsingular transformation of the Polish probability space $(Y,\nu)$ and $\a\subset\B(Y)$ is a countable partition so that for each $a\in\a,\ \tau a\in\s(\a)$ and $\tau:\a\to\tau a$ is nonsingular, invertible.
\

It follows that if $(Y,\nu,\tau,\a)$  is a Markov map, then so is
$(Y,\nu,\tau^n,\a_n)\ \forall\ n\ge 1$ where $\a_n:=\bigvee_{k=0}^{n-1}\tau^{-k}\a$.
\

The transfer operators $\tau^n:L^1(\nu)\hookleftarrow$ ($n\ge 1$) are given by
$$\tau^n f=\sum_{a\in\a_n}1_{\tau^n a}v_a'f\circ v_a$$
where for $a\in\a_n,\ v_a:\tau^n a\to a,\ v_a\circ\tau=\text{\tt Id}$ on $a$ and
$v_a':=\tfrac{d\nu\circ v_a}{d\,\nu}$.
\

As in  \cite{AD} (also \cite{A-Sera}), a  {\it Gibbs-Markov map} is a  {\it Markov map} 
  $(Y,\nu,\tau,\a)$ so that
  $\inf_{a\in\alpha}\nu(\tau a)>0$ 
 and, for some $\th\in (0,1)$
 \begin{align*}\sup_{n\ge 1,\ a\in\alpha_n,\ x,y\in \tau^na}
 \frac1{\th^{t(x,y)}}\cdot \left|\frac{v_a'(x)}{v_a'(y)}-1\right|<\infty
\end{align*}
where $t(x,y)=\min\{n\ge 1:\ \a_n(x)\ne \a_n(y)\}\le\infty$ with $x\in\a_n(x)\in\a_n$
\

If $(Y,\nu,\tau,\a)$  is a Gibbs-Markov map, then so is $(Y,\nu,\tau^n,\a_n)\ \forall\ n\ge 1$.
\

As in  \cite{AD} (also \cite{A-Sera}),  the transfer operator $\ttau$ acts quasicompactly on the space of   $(\a,\th)$-{\it H\"older}\label{Holder} functions on $Y$; that is
$$L_{\a,\th}:=\{f:Y\to\Bbb R:\ D_{Y,\th,\a}(f):=\sup_{x,y\in \Om}\tfrac{|f(x)-f(y)|}{\th^{t(x,y)}}<\infty\}$$ 
equipped with the norm $\|f\|_{L_{\a,\th}}:=\|f\|_1+D_{Y,\th,\a}(f)$.

\

Let  $(X,m,T)$ be a  conservative, ergodic,  measure preserving 
transformation . We'll call the set $\Om\in\B(X),\ 0<m(\Om)<\infty$  a {\it Gibbs-Markov set} if there is a countable partition $\a\subset\B(\Om)$ so that 
$(\Om,m_\Om,T_\Om,\a)$ is a mixing Gibbs-Markov map  and the first return time $\v_\Om:\Om\to\Bbb N$ is $\a$-measurable. In this case, a $(n,\a)$-{\it cylinder} ($n\ge 1$) is a set of form
$$a=\bigcap_{j=0}^{n-1}T_\Om^{-j}a_j$$ where  $a_0,\dots,a_{n-1}\in\a$. The collection of $(n,\a)$-cylinders  is $\a_n:=\bigvee_{j=0}^{n-1}T_\Om^{-j}\a$ and the collection of cylinders is 
$$\mathcal{C}_\a:=\bigcup_{n=1}^\infty\a_n.$$
\

\subsection*{Ergodic functional tied-down renewal properties}
\

These (i.e. \ref{umbrella} and \ref{knight} below) are functional analogues of the tied-down renewal mixing properties considered in \cite{A-Sera}.
\

   If the conservative, ergodic,  measure preserving 
transformation  $(X,m,T)$ has a Gibbs-Markov set $\Om$, then it is pointwise dual ergodic, whence rationally ergodic (with the same return sequence); and weakly mixing iff $\v_\Om$ is non-arithmetic with respect to $T_\Om$. It is standard that these latter properties lift to the natural extension.
    \
    
\

By \cite[Proposition 1.1]{RatErg}, if $(X,m,T)$ is a rationally ergodic measure preserving
transformation , then there is a dense, $T$-invariant, hereditary ring $R(T)\subset \mathcal F_+$ (as on p. \pageref{mathcalF}) and a sequence $a(n)=a_n(T)>0$ (the {\it return sequence}) so that
\begin{align*}\tag*{{\scriptsize \faTree}}\label{faTree}\tfrac1{a(n)}\sum_{k=0}^{n-1}m(A\cap T^{-k}B)
\xrightarrow[n\to\infty]{×}m(A)m(B)\ \forall\ A,B\in R(T). 
\end{align*}
A measure preserving
transformation  $(X,m,T)$ satisfying \ref{faTree} is called (in \cite{RatErg}) {\it weakly rationally ergodic}.

\

Our main results are integrated, functional versions of \cite[Theorem A]{A-Sera}.
\

In the following,   $\Psi_{n,f}(t):=\frac{S_{\lfl nt\rfl}(f)}{a(n)}$
($f\in L^1_+$).
 \proclaim{Theorem U}\label{thmU}
\ 

\bul Suppose that $(X,m,T)$ is a  pointwise dual ergodic  conservative, ergodic,  measure preserving 
transformation  with $a(n)=a_n(T)\ \g$-regularly 
 varying with $0<\g<1$, then

 with $u(n):=\tfrac{\g a(n)}n$,
 \begin{align*}\tag*{{\small\faUmbrella}}\label{umbrella} 
\tfrac1{a(n)}\sum_{k=1}^n&\int_Ag(\Psi_{k,f})1_B\circ T^kdm\xrightarrow[n\to\infty]{}
m(A)m(B)\Bbb E(g(m(f)\frak w_\g))\\ & \forall\ A,\ B\in R(T),\ f\in L^1_+,\ g\in C_{B}( D_{\uparrow,1}).\end{align*}
\bul\ \ If, in addition, $(X,m,T)$ is weakly mixing and   has a Gibbs-Markov set $\Om\in\mathcal{F}_+$, then
\begin{align*}\tag*{\symknight}\label{knight}
\tfrac1{a(n)}\sum_{k=1}^n&|\T^{k}(1_Ag(\Psi_{k,f}))-u(k)m(A)\Bbb E(g(m(f)\mathfrak w_\g))|\\ & \xrightarrow[n\to\infty]{}\ 0\ \ \ \text{a.e.}\ 
\forall\ A\in\mathcal F_+,\ f\in L^1_+,\ g\in C_{B}( D_{\uparrow,1}).\end{align*}
\bul\  If, in addition $(\Om,m_\Om,T_\Om,\a)$ satisfies
 \begin{align*}\tag{{\tt OSRT}}\label{OSRT}
	 \tfrac1{u(n)}\T^{n}(1_A)
	 \xrightarrow[n\to\infty]{}\  m(A)\ \text{uniformly on}\ \Om\  \forall\ A\in\mathcal C_\a,
\end{align*}
then
\begin{align*}\tag*{\Industry}\label{Industry}
	  \tfrac1{u(n)}\T^{n}(1_Ag(&\Psi_{n,1_\Om}))
	 \xrightarrow[n\to\infty]{} \mathbb{E}(g(\mathfrak w_\g))m(A)\ \ \text{uniformly on}\ \Om\\ &   \forall\ A\in\mathcal C_\a, 
	 g\in C_B(D_{\uparrow,1}).
\end{align*}
\endproclaim
\f{\bf Remarks}\ \ 
\

\sms {(i)}\ It is standard to show that if a non-invertible conservative, ergodic,  measure preserving 
transformation  satisfies \ref{umbrella}, then so does its natural extension. Also, if $T$ satisfies \ref{knight}, then $T$ and its natural extension also satisfy
\begin{align*}
\tfrac1{a(n)}\sum_{k=1}^n&|\int_{A\cap T^{-k}B}g(\Psi_{k,f})dm-u(k)m(A)m(B)\Bbb E(g(m(f)\mathfrak w_\g))|\\ & \xrightarrow[n\to\infty]{}\ 0\ \ \ \
\forall\ A,\ B\in R(T),\ f\in L^1_+,\ g\in C_{B}( D_{\uparrow,1}).\end{align*}
\endproclaim
\

\sms (ii)
The condition (\ref{OSRT}) for  $(\Om,m_\Om,T_\Om,\a)$ was shown in \cite{MelTer} for $\tfrac12<\g<1$ (with no further assumptions) and in \cite{Goue11}
  for $0<\g<1$  under the additional assumption that the $\a$-measurable return time function $\v:\Om\to\Bbb N$ is aperiodic and 
 $m_\Om([\v=n])\ll \tfrac{m_\Om([\v\ge n])}n$. 
 \

\subsection*{Examples: Intermittent interval maps}
 \

An {\tt\small intermittent interval map} satisfying  {\tt\small Thaler's conditions} as in \cite{Tha83}:
\sbul admits an absolutely continuous invariant
measure with density continuous away from the indifferent fixed points (\cite{Tha80});
\sbul admits {\tt\small wandering rates} (\cite[\S3]{Tha83});
\sbul is pointwise dual ergodic and 
 has Gibbs-Markov  sets (\cite[Theorem 3]{RFEx}). See also \cite[Chapter 4]{IET}.
\

In particular, each map $T_\g:[0,1]\to [0,1]$ ($0<\gamma<1$) defined  by
\begin{align*}
	T_\g x
	=
	\begin{cases}
		x(1+(2x)^{1/\gamma}), &0\leq x<\frac{1}{2},
		\\
		2x-1, &\frac{1}{2}\leq x\leq 1
	\end{cases}
\end{align*}
 satisfies Thaler's conditions as above and indeed all  the assumptions of Theorem U.\footnote{and are aka
 ``{\tt LSV maps}'' having been  considered in \cite{LSV} for the different  parameters $\g>1$.}
 \
 
 Thus if  $h_\g$ is the $T_\g$-invariant density,  then $([0,1],\mu^\g, T_\g)$ is pointwise dual ergodic
 with  $d\mu^\g(x):=h_\g(x)dx$.
 \
 
 It follows from the regularly varying expansion near 
the indifferent fixed point $0$ that the invariant density $h_\g$ for $T_\g$ 
 satisfies $h_\g(x)\propto  x^{-1/\gamma}$ as $x\to 0$, whence the  wandering rate $L_{T_\g}(n)\propto n^{1-\g}$ and the return sequence $a_n(T_\g)\propto n^{\gamma}$. Thus  
 the functional Darling-Kac theorem holds.
\

The set 
 $\Omega=[\frac{1}{2},1]$ is a Gibbs-Markov set for each $T_{\g}$ with respect to the return time partitions $\alpha_\g:=\{[\v_\g=n]\cap\Om:\ n\ge 1\}$ (where
 $\v_\g(x):=\min\{j\ge 1:\ T_\g^j(x)\in\Om\}$.  Moreover $\mu^\g_\Om([\v_\g=n])\propto\tfrac1{n^{1+\g}}$
whence $(\Omega, \mu^\g_\Omega, (T_\g)_\Omega, \alpha)$ satisfies (\ref{OSRT}).
Therefore each $([0,1], \mu^\g, T_\g)$ satisfies all the assumptions of Theorem U.
\subsection*{Outline of Proofs}
\

The proofs of Propositions C, D, E (p.\pageref{pfCD}) and \ref{umbrella} (p.\pageref{umbrella}) in Theorem U, are interwined and use
\bul Proposition 4 (p.\pageref{prop4}) which ``converts'' a distibutional limit to a tied down  limit; and 

\bul Lemma 5 (p.\pageref{lem5}) which gives tied down limits for well behaved renewal shifts.

 The  statement \ \ref{knight}\  in Theorem U will follow (in \S3) from
 the Uniform GL Lemma ((\ref{UGL}) on p. \pageref{UGL}). \footnote{ "GL" stands for Garsia-Lamperti in honor of \cite[Theorem 1.1]{G-L} which inspired (UGL).}
\section*{\S2 Tied down Mittag-Leffler limits}

\proclaim{Proposition 4}\label{prop4} \ Let $(X,m,T)$ be a weakly rationally ergodic measure preserving
transformation  \ \ \footnote{as in \ref{faTree} on p. \pageref{faTree}} with $\g$-regularly varying return sequence $a(n)=a_n(T)$  ($\g\in (0,1]$).
\

Suppose that \begin{align*}\tag*{\symrook}\label{rook}\Psi_{n,f}
\xrightarrow[n\to\infty]{\frak d}m(f)\text{\tt m}_\g\ \text{\rm in}\ (D_{\uparrow,\infty},J_1)\ \forall\ f\in L^1(m)_+\end{align*}
where
$\Psi_{n,f}(t):=\frac{S_{\lfl nt\rfl}(f)}{a(n)}$ and  $\text{\tt m}_\g\in\text{\tt RV}\,( D_{\uparrow,\infty})$ 
is a.s.  continuous with $\text{\tt m}_\g(0)=0$, then 

\begin{align*}\tag*{ \sympawn}\label{pawn} 
\tfrac1{a(n)}\sum_{k=1}^n&\int_Ag(\Psi_{k,f})1_\Om\circ T^kdm\xrightarrow[n\to\infty]{}
m(A)m(\Om)\Bbb E(g(m(f)\text{\tt w}(\text{\tt m}_\g)))\\ & \forall\ A,\ \Om\in R(T),\ f\in L^1_+,\ g\in C_{B,c}( D_{\uparrow,1})\end{align*}
where $\text{\tt w}(\text{\tt m}_\g)\in\text{\tt RV}\,( D_{\uparrow,1})$ 
satisfies:
\begin{align*}\tag*{{\scriptsize\faBook}}\label{bcbook} \Bbb E(h(\text{\tt w}(\text{\tt m}_\g)))=\Bbb E(\int_0^1h(\D_{t,\g}\text{\tt m}_\g)d\text{\tt m}_\g(t))\ \forall\ h\in C_{B}( D_{\uparrow,1}).\end{align*}
\endproclaim

\demo{Proof}\ \ We assume WLOG that $g=\mathbb{g}_{\u t}\in\Pi$ (as in \ref{dsliterary} on p. \pageref{dsliterary}) with  $x\mapsto \log h(e^x)\ \&\ x\mapsto \log g_\nu(e^x)$ ($1\le\nu\le N$)  uniformly continuous.

\Par1\ \ We first show \
\ref{pawn} \ with $A=\Om\in R(T)\ \&\ f=1_\Om$. For convenience we assume $m(\Om)=1$. 
Fix $0<\mathcal E<1$. 
\

Writing $s_n:=S_n(1_\Om)\ \&\ \psi_n(t):=\Psi_{n,1_\Om}(t)=\tfrac{s_{\lfl nt\rfl}}{a(n)}$, we have
 uniformly in $t\in [\mathcal E,1]$, with 
$k=\lfl nt\rfl\ \&\ u\in (0,1)$, that 
\begin{align*}\psi_k(u)&=\tfrac{s_{ku}}{a(k)}=\tfrac{s_{\lfl\lfl nt\rfl u\rfl}}{a(nt)}\\ &\sim \tfrac{s_{\lfl ntu\rfl}}{t^\g a(n)}
\\ &=
\tfrac{\psi_n(tu)}{t^\g}=:(\D_{t,\g}\psi_n)(u).
\end{align*}
            
Next,
\begin{align*}\tag*{\dsagricultural}\label{dsagricultural}\int_\mathcal E^1g(\D_{t,\g}\psi_n)d\psi_n(t)&\sim\sum_{\frac{k}n\in (\mathcal E,1)}
\int_{\frac{k}n}^{\frac{k+1}n}g(\D_{t,\g}\psi_n)d\psi_n(t)\\ &\sim
\tfrac1{a(n)}\sum_{\frac{k}n\in (\mathcal E,1)}g(\psi_k)1_\Om\circ T^k.
\end{align*}

For $\xi\in  D_{\uparrow,1},\ 0<\mathcal E<1$ 
, let
$$\frak L_{g,\mathcal E}(\xi):=\int_\mathcal{E}^1g(\D_{t,\g}\xi)d\xi(t).$$
We claim that
\sms\  \ If  $\xi_n,\ \xi\in   D_{\uparrow,1}$ with $\xi$ continuous and  $\xi_n\xrightarrow[n\to\infty]{\text{\tt\tiny DF}}\ \xi$, then for
$0<\mathcal E<1$,
\begin{align*}\tag*{{\scriptsize\faCogs}}\label{faCogs}
 \frak L_{g,\mathcal E}(\xi_n)\xrightarrow[n\to\infty]{}\ \frak L_{g,\mathcal E}(\xi).
\end{align*}

\demo{Proof of \ref{faCogs} }\ Fix $0<\mathcal E<1$.
\

Since $\xi$ is continuous, so is each $\D_{t,\g}\xi$ ($t\in (0,1]$) whence  the mapping $t\mapsto \D_{t,\g}\xi$ is continuous $[\mathcal{E},1]\to (D_{\uparrow,1},J_1)$.
\

The collection $\{\D_{t,\g}\xi:\ t\in [\mathcal{E},1]\}$ is therefore compact in
$(D_{\uparrow,1},J_1)$  and for each $\e>0,\ \exists\ \d(\e)>0$ so that
$$\zeta\in D_{\uparrow,1},\ t\in [\mathcal{E},1]:\ \sup_{s\in [0,1]}|\zeta(s)-\D_{t,\g}\xi(s)|<\d(\e)\  \Rightarrow\ |g(\zeta)-g(\D_{t,\g}\xi)|<\e.$$
\

Thus, if  $\zeta\in D_{\uparrow,1}$ and   $\sup_{s\in [0,1]}|\zeta(s)-\xi(s)|<\d(\e)\mathcal{E}^\g$, then for each $t\in [\mathcal{E},1],\ s\in [0,1]$, we have
\begin{align*}|\D_{t,\g}\zeta(s)-\D_{t,\g}\xi(s)|\le\tfrac1{\mathcal{E}^\g}|\zeta(st)-\xi(st)|<
\d(\e) 
\end{align*}
whence
$$\sup_{t\in [\mathcal{E},1]}|g(\D_{t,\g}\zeta)-g(\D_{t,\g}\xi)|<\e.$$
Thus, since $\xi_n\xrightarrow[n\to\infty]{×}\xi$ uniformly on $[0,1]$,
$$\d(\mathcal E,n):=\sup_{t\in [\mathcal E,1]}|g(\D_{t,\g}\xi_n)-g(\D_{t,\g}\xi)|\xrightarrow[n\to\infty]{×}0$$
and
\begin{align*}
 \int_{\mathcal{E}}^1g(\D_{t,\g}\xi_n)d\xi_n(t)=\int_{\mathcal{E}}^1g(\D_{t,\g}\xi)d\xi_n(t)\pm \d(\mathcal{E},n)\xi_n(1).
\end{align*}
Next, since $t\mapsto g(\D_{t,\g}\xi)$ is continuous $[\mathcal E,1]\to \Bbb R$,
$$\int_{\mathcal{E}}^1g(\D_{t,\g}\xi)d\xi_n(t)\xrightarrow[n\to\infty]{}\int_{\mathcal{E}}^1g(\D_{t,\g}\xi)d\xi(t).$$
Thus
$$\frak L_{g,\mathcal E}(\xi_n)\xrightarrow[n\to\infty]{}\frak L_{g,\mathcal E}(\xi).\ \ \CheckedBox\ \text{\ref{faCogs}}$$
\

 In view of \ref{faCogs},  we now have by Skorokhod's Representation Theorem 
(\cite[Thm. 6.7]{Billingsley})  and \ref{dsagricultural} that $\forall\ 0<\mathcal{E}<1$, in 
$\text{\tt RV}(D_{\uparrow,1},J_1)$,
$$\tfrac1{a(n)}\sum_{\frac{k}n\in (\mathcal E,1)}g(\psi_k)1_\Om\circ T^k\sim\frak L_{g,\mathcal E}(\psi_n)\xrightarrow[n\to\infty]{}\frak L_{g,\mathcal E}(\text{\tt m}_\g).$$
To complete the proof of \P1, we must establish the limits as $\mathcal{E}\to 0$.
\

To this end we claim that $\Bbb E(\text{\tt m}_\g(1))\le 1$. Indeed
\begin{align*}
 1&\xleftarrow[n\to\infty]{}\ 
 \int_\Om\tfrac{s_n}{a(n)}dm\ \because\ \Om\in R(T),\ m(\Om)=1 \\ &\ge \int_\Om(\tfrac{s_n}{a(n)}\wedge M)dm\ \forall\ M>0,\\ &
 \xrightarrow[n\to\infty]{}\Bbb E(\text{\tt m}_\g(1)\wedge M)\xrightarrow[M\uparrow\infty]{}
 \Bbb E(\text{\tt m}_\g(1)).
\end{align*}
In view of this, since $\text{\tt m}_\g$ is a.s. continuous $\&\ \text{\tt m}_\g(0)=0$, 
by dominated convergence
$$\Bbb E(\text{\tt m}_\g(\mathcal{E}))\xrightarrow[\mathcal{E}\to 0+]{×}0,$$
whence
$$\Bbb E(\int_0^\mathcal{E}g(\D_{t,\g}\text{\tt m}_\g)d\text{\tt m}_\g)\le \|g\|_{C_B} \Bbb E(\text{\tt m}_\g(\mathcal{E}))\xrightarrow[\mathcal{E}\to 0+]{×}0.$$
Moreover, by the $\g$-regular variation of $a(n)$,
$$\tfrac1{a(n)}\sum_{\frac{k}n\in (0,\mathcal E)}\int_\Om g(\psi_k)1_\Om\circ T^kdm\lesssim
\|g\|_{C_B}\tfrac{a(n\mathcal E)}{a(n)}\sim\ \mathcal E^\g\|g\|_{C_B}.$$\footnote{Here and throughout, for $a_n,\ b_n>0,\ 
a_n\lesssim b_n$ means $\varlimsup_{n\to\infty}\frac{a_n}{b_n}\le 1$.}

\
This proves  \ref{pawn} with $A=\Om\in R(T)\ \&\ f=1_\Om$. \ \Checkedbox\ \P1

\

\Par2\ \ Next, we prove  \ref{pawn} for any $A,\ \Om\in R(T)$ with $f=1_\Om$.
\

Define $G_N:\Om\to\Bbb R_+$ by
$$G_N:=\tfrac1{a(N)}\sum_{k=1}^Ng(\psi_k)1_\Om\circ T^k,$$
then
\bul $\int_\Om G_Ndm\xrightarrow[N\to\infty]{}\Bbb E(g(\text{\tt w}(\text{\tt m}_\g)))$ by \P1;
\bul $G_N\le \|g\|_{C_{B,c}}\tfrac{s_N}{a(N)}$, whence $\{G_N:\ N\ge 1\}$ is weakly, sequentially
compact in $L^1(F)\ \forall\ F\in R(T)$; and
\bul $G_N-G_N\circ T\xrightarrow[N\to\infty]{m}\ 0$.
\

Fix $F\in R(T),\ F\supset\Om\cup A$. Suppose that $G_{N_k}\xrightarrow[k\to\infty]{}\ H$ weakly in $L^1(F)$,
then $H=H\circ T$ whence by ergodicity, $H$ is constant and
$$H=\int_\Om Hdm\xleftarrow[k\to\infty]{}\int_\Om G_{N_k}dm\xrightarrow[k\to\infty]{}\Bbb E(g(\text{\tt w}(\text{\tt m}_\g))).$$
Thus
$$G_N\xrightarrow[N\to\infty]{}\ \Bbb E(g(\text{\tt w}(\text{\tt m}_\g)))\ \text{weakly in }\ L^1(F)$$
which establishes \P2.

\

\Par3\ \ To finish, we show \ref{pawn} for any $A,\ \Om\in R(T)$ and
$f\in L^1(m)_+$.
\

\demo{Proof of \P3}\label{pfP3} \ \ 
\

Fix $A,\ \Om\in R(T)\ \&\ f\in L^1(m)_+$.

\

By the ratio  ergodic theorem and Egorov's theorem, 

$\exists\ A_\nu\in\B(A),\ \ A_\nu\uparrow\ A$\  so that for each $\nu\ge 1$,
$$\tfrac{S_n(f)}{s_n}\xrightarrow[n\to\infty]{} \tfrac{m(f)}{m(\Om)}\ \text{uniformly on}\ A_\nu$$
whence 
$$\sum_{k=1}^n1_\Om\circ T^kg(\Psi_{n,f})\sim\ \sum_{k=1}^n1_\Om\circ T^kg( \tfrac{m(f)}{m(\Om)}\psi_n)\ \text{uniformly on}\ A_\nu.$$
Thus
\begin{align*}\tfrac1{a(n)}&\sum_{k=1}^n\int_{A_\nu\cap T^{-k}\Om}g(\Psi_{n,f})dm\\ &\approx\ \tfrac1{a(n)}\sum_{k=1}^n\int_{A_\nu\cap T^{-k}\Om}g( \tfrac{m(f)}{m(\Om)}\psi_{n})dm\\ &\xrightarrow[n\to\infty]{}
m(A_\nu)m(\Om)\Bbb E(g(m(f)\text{\tt w}(\text{\tt m}_\g))) 
\end{align*}\footnote{where $a_n\approx b_n$ means $a_n-b_n\xrightarrow[n\to\infty]{×}0$}
\

by \P2  and
\begin{align*}\big|\tfrac1{a(n)}&\sum_{k=1}^n\int_{A\cap T^{-k}\Om}g(\Psi_{n,f})dm-m(A)m(\Om)\Bbb E(g(m(f)\text{\tt w}(\text{\tt m}_\g)))\big|
 \\ &\lesssim\ \|g\|_{C_B}(\tfrac1{a(n)}\sum_{k=1}^nm(A\setminus A_\nu\cap T^{-k}\Om)+m(A\setminus A_\nu)m(\Om))\\ &\xrightarrow[n\to\infty]{}\ 2\|g\|_{C_B}m(A\setminus A_\nu)m(\Om)\\ &
 \xrightarrow[\nu\to\infty]{}\ 0.\ \CheckedBox\ \P3
\end{align*}
This completes the proof of \ref{pawn}  and Proposition 4.\ \ \CheckedBox

\proclaim{Lemma 5}\label{lem5}
\

\ \ Let $(X,m,T)=(\Om,p,S)^\phi$ be a renewal shift with aperiodic lifetime distribution
 $f\in\mathcal P(\Bbb N)$ so that
$a(n)$ is $\g$-regularly varying with $\g\in (0,1)$ and
satisfying (\ref{SRT}) on p. \pageref{SRT}, 
then

 for $\mathbb{g}_{\u t}\in\Pi$ a product function
as in \ref{dsliterary} {\rm (p. \pageref{dsliterary})},
\begin{align*}\tag*{\symking}\label{king}
 \tfrac1{u(N)}\int_{\widetilde{\Om}\cap T^{-N}\widetilde{\Om}} \mathbb{g}_{\u t}(\psi_{N})dm\xrightarrow[N\to\infty]{}\frak e_\g(\mathbb{g}_{\u t}).
\end{align*}
where  $\frak e_\g(\mathbb{g}_{\u t})$ denotes the RHS of \ref{queen} on p. \pageref{queen}.\endproclaim
\subsection*{Continuity of functionals}
\

The proofs  below  need  \cite[Section 5]{Toru-thesis}:

\

For  $\Phi=\mathcal{I},\ \mathcal{G}\circ\mathcal{I},\ \mathcal{D}\circ\mathcal{I}$ (as on p. \pageref{waiting}),
\begin{align*}\tag*{{\Large\smiley}}\label{Smi}  \text{If}\ \xi_n\xrightarrow[n\to\infty]{J_1}&\xi,\ 
\xi(0)=0,\ \xi\   \text{strictly increasing}\ \&\ \xi(\infty-)=\infty,\\ &  \text{then}\ 
    \Phi(\xi_n)\xrightarrow[n\to\infty]{J_1}\Phi(\xi). 
\end{align*}

\subsection*{Proof of Lemma 5}\label{proofP5}
\

Let
\begin{align*}\tag*{\faThumbTack}\label{faThumbTack}
\Pi_0:=\{\mathbb{g}_{\u t}\in\Pi:\ N\geq1,\ \u t\in\Bbb{Q}^N_+,\ g_1,g_2,\dots,g_N,h\in \mathcal{A}\}
\end{align*}
with $\mathbb{g}_{\u t}\ \&\ \Pi$  as in \ref{dsliterary} (p.\pageref{dsliterary}), and where $\mathcal{A}\subset C([0,\infty],[0,1])$ is countable, dense, $\mathbb{1}\in\mathcal{A}$ and  $g\in\mathcal{A}\Rightarrow 1-g\in\mathcal{A}$.

We'll first establish that 
\begin{align*}\tag*{\faFlash}\label{faFlash}
 \varliminf_{N\to\infty}\tfrac1{u(N)}\int_{{\widetilde{\Om}}\cap T^{-N}{\widetilde{\Om}}} \mathbb{g}_{\u t}(\psi_{N})dm
 \ge\frak e_\g(\mathbb{g}_{\u t})\ \forall\ \mathbb{g}_{\u t}\in\Pi_0.
\end{align*}
To this end, define $Z_n(x):=\min\,\{k\ge n+1:\,T^kx\in{\widetilde{\Om}}\}$, then   
$$\mathcal D(\psi_n)(t)=\tfrac1nZ_{\lfl nt\rfl}.$$
By \cite{TK-cgce},\ $\psi_n^{-1}\xrightarrow[n\to\infty]{\text{\tt\tiny RV}(D_{\uparrow,\infty},J_1)}\frak m_\g^{-1}$, whence by \ref{Smi} on p. \pageref{Smi}:
\begin{align*}\tag*{\dscommercial}\label{dscommercial}
(\psi_n,\mathcal D(\psi_n)(t_N))\xrightarrow[n\to\infty]{\text{\tt\tiny RV}((D_{\uparrow,\infty},J_1)\x\Bbb R_+)}(\frak m_\g,\mathcal D(\frak m_\g)(t_N))
\end{align*}

Next, writing $\kappa:=\lfl nt_N\rfl$,
\begin{align*}1_{\widetilde{\Om}}\mathbb{g}_{\u t}(\psi_n)1_{\widetilde{\Om}}  \circ T^n&=
\sum_{k=\kappa+1}^n1_{\widetilde{\Om}}\mathbb{g}_{\u t}(\psi_n)1_{[Z_\kappa=k]}1_{\widetilde{\Om}}  \circ T^n\\ &=
\sum_{k=\kappa+1}^n1_{\widetilde{\Om}}\frak g(\psi_n)1_{[Z_\kappa=k]}1_{\widetilde{\Om}}\circ T^k h(\tfrac{s_{n-k}\circ T^k}{a(n)})1_{\widetilde{\Om}}  \circ T^n.
\end{align*}

Fix $\l\in (t_N,1)$. Using the Markov property, we have 

\begin{align*}\int_{\widetilde{\Om}\cap T^{-n}\widetilde{\Om}}\mathbb{g}_{\u t}(\psi_n)dm&=
\sum_{k=\kappa+1}^n\int_{\widetilde{\Om}} \frak g(\psi_n)1_{[Z_\kappa=k]}dm\int_{\widetilde{\Om}} h(\tfrac{s_{n-k}\circ T^k}{a(n)})1_{\widetilde{\Om}}  \circ T^{n-k}dm\\ &\ge 
\sum_{\kappa+1\le k\le n\l }
\\ &=:I_n(\l).
\end{align*}
To establish \ref{faFlash}, we'll show that $\exists\ L:(t_N,1)\to\Bbb R_+$ so that
\begin{align*}&\tag*{{\Large\Chair}}\label{Chair} \tfrac{I_n(\l)}{u(n)}\xrightarrow[n\to\infty]{×}L(\l)\xrightarrow[\l\to 1-]{×}\frak e_\g(\mathbb{g}_{\u t}).
\end{align*}
\demo{Proof of \Chair}
\

By  (\ref{SRT}) 
on p. \pageref{SRT}, $m(\widetilde{\Om}\cap T^{-n}\widetilde{\Om})= u(n)\sim\tfrac{\g a(n)}n$, whence by \cite[lemma 2.1]{A-Sera} and the ``remarks about mixing'' there, \ \ 
\begin{align*}\tag*{\faRocket}\label{faRocket}
\int_{\widetilde{\Om}\cap T^{-n}\widetilde{\Om}} h(\tfrac{S_n(1_\Om)}{a(n)})dm= \int_{\widetilde{\Om}\cap T^{-n}\widetilde{\Om}} h(\psi_{n}(1))dm\sim \Bbb E(h( W_\g ))\cdot u(n).
\end{align*}

For $\rho>1$, fix  $K=K_\rho\ge 1$ so that
$$\int_{\widetilde{\Om}\cap T^{-N}\widetilde{\Om}}h(\psi_{N}(1))dm=\rho^{\pm 1} \Bbb E(h( W_\g ))\cdot u(N)\ \forall\ N\ge K,$$
then for fixed $\l\in (t_N,1)\ \&\ n\ge \tfrac{K}{1-\l}$,
\begin{align*}I_n(\l)&=\rho^{\pm 1}
\sum_{\kappa+1\le k\le n\l }\int_{\widetilde{\Om}} \frak g(\psi_n)1_{[Z_\kappa=k]}dm\cdot \Bbb E(h((1-\tfrac{k}n)^\g W_\g))u(n-k)\\ &=
\rho^{\pm 1}\sum_{\kappa+1\le k\le n\l }\int_{\widetilde{\Om}} \frak g(\psi_n)1_{[Z_\kappa=k]}dm\cdot (\tfrac1{1-\frac{k}n})^{1-\g}\Bbb E(h((1-\tfrac{k}n)^\g W_\g))u(n)\ \ \text{\small for large $n$ with increased $\rho$ by regular variation of}\ u\\ &=
\rho^{\pm 1}\int_{\widetilde{\Om}} \frak g(\psi_n)1_{[Z_\kappa\le n\l ]}(\tfrac1{1-\mathcal D(\psi_n)(t_N)})^{1-\g} \Bbb E(h((1-\mathcal D(\psi_n)(t_N))^\g W_\g))dm\cdot u(n).
\end{align*}
By \ref{dscommercial} on p. \pageref{dscommercial},
\begin{align*}
\int_{\widetilde{\Om}} &\frak g(\psi_n)1_{[Z_\kappa\le n\l ]}(\tfrac1{1-\mathcal D(\psi_n)(t_N)})^{1-\g}\cdot \Bbb E(h((1-\mathcal D(\psi_n)(t_N))^\g W_\g))dm\\ &\xrightarrow[n\to\infty]{×}\Bbb E\(\frak g(\frak m_\g)1_{[\mathcal D(\frak m_\g)(t_N)\le \l]}(\tfrac1{1-\mathcal D(\frak m_\g)(t_N)})^{1-\g}h((1-\mathcal D(\frak m_\g)(t_N))^\g W_\g) \)
\\ &=:L(\l)\xrightarrow[\l\to 1-]{×}\frak e_\g(\mathbb{g}_{\u t})\ \ \text{establishing \ref{Chair}.\ \ \CheckedBox}.
\end{align*}
\demo{Proof of}\ \ref{king}
\

Evidently  \ref{faFlash} (on p. \pageref{faFlash})  holds  $\forall\ g\in\Pi_1$ where
\begin{align*}\tag*{\faLock}\label{faLock}
\Pi_1:=\{\sum_{m=1}^M \mathbb{g}^{(m)}:\ M\geq1,\  \mathbb{g}^{(1)},\mathbb{g}^{(2)},\dots\mathbb{g}^{(M)}\in \Pi_0\}	
\end{align*}
with $\Pi_0$ as in \ref{faThumbTack} (p.\pageref{faThumbTack})
and where $\frak e_\g(\sum_{m=1}^M \mathbb{g}^{(m)}):=\sum_{m=1}^M \frak e_\g(\mathbb{g}^{(m)})$.
\

Moreover, $\mathbb{g}_{\u t}\in \Pi_0 \Rightarrow 1-\mathbb{g}_{\u t}\in \Pi_1$ because, for $\mathbb{g}_{\u t}\in \Pi_0$ and $\xi\in D_{\uparrow,1}$,
\begin{align*}
  1-\mathbb{g}_{\u t}(\xi)=1-\prod_{j=1}^Ng_j\cdot h=
  \sum_{k=1}^N(1-g_k)\prod_{k<j\le N}g_jh+1-h
\end{align*}
with each $g_k:=g_k(\xi(t_k)),\ 1-g_k,\  h:=h(\xi(1)-\xi(t_N)),\ 1-h\in\mathcal{A} .$

For $\mathbb{g}_{\u t}\in\Pi_0$, we have on ${\widetilde{\Om}}$:

\begin{align*}
\frak e_\g(\mathbb{g}_{\u t}) &\lesssim \tfrac1{u(n)}\T^{n}(1_{\widetilde{\Om}}\mathbb{g}_{\u t}(\psi_{n})) \ \text{by\ \ref{faFlash}}\\ &
=\tfrac1{u(n)}\T^{n}(1_{\widetilde{\Om}})-\tfrac1{u(n)}\T^{n}(1_{\widetilde{\Om}}(1-\mathbb{g}_{\u t}(\psi_{n})))\\ &
=
1-\tfrac1{u(n)}\T^{n}(1_{\widetilde{\Om}}(1-\mathbb{g}_{\u t}(\psi_{n}))) \ \because\ \T^n1_{\widetilde{\Om}}=u(n)\ \text{on}\ {\widetilde{\Om}}
 \\
 &\lesssim 1-\frak e_\g(1-\mathbb{g}_{\u t}) \ \text{by\  \ref{faFlash}}
 \\
 &=\frak e_\g(\mathbb{g}_{\u t}).
\end{align*}

This proves \ref{king} and Lemma 5.\ \ \Checkedbox

\

Proposition E follows from Lemma 5 and Proposition D.
\demo{Proof of Propositions C and D}\label{pfCD}
\

Fix $\g\in (0,1)$.  Let, as in Lemma 5, $(X,m,T)=(\Om,p,S)^\phi$ be a renewal shift with aperiodic lifetime distribution
 $f\in\mathcal P(\Bbb N)$ so that
$a(n)$ is $\g$-regularly varying with $\g\in (0,1)$ and which satisfies (\ref{SRT}).

Let
$$\psi_n(t):=\tfrac{s_{[nt]}}{a(n)},\ s_n:=S_n(1_{\widetilde{\Om}}),\   \&\ z_n:=\max\{k\le n:\ T^kx\in{\widetilde{\Om}}\},$$ then
$\tfrac{z_n}n=\mathcal G(\psi_n)(1)$.

By \cite{TK-cgce,TK-weak},  (see \ref{faMedkit} on p. \pageref{faMedkit}), 
$$\psi_n^{-1}\xrightarrow[n\to\infty]{\text{\tiny\tt RV}(D_{\uparrow,\infty},J_1)}\eta_\g=\frak m_\g^{-1}.$$
Thus  by \ref{Smi} on p. \pageref{Smi}, $(\psi_n,\frac{z_n}n)\xrightarrow[n\to\infty]{\text{\tt\tiny RV}((D_{\uparrow,\infty},J_1)\x\Bbb R_+)}(\frak m_\g,\mathcal G(\frak m_\g)(1))$ whence 
$$(\D_{\frac{z_n}n,\g}\psi_n)1_{[\frac{z_n}n\ne 0]}\xrightarrow[n\to\infty]{\text{\tiny\tt RV}(D_{\uparrow,1},J_1)}\D_{\mathcal G(\frak m_\g)(1),\g}\frak m_\g=:\frak w_\g$$
and for  for $\mathbb{g}_{\u t}\in\Pi_0$ 
as in \ref{faThumbTack} on p. \pageref{faThumbTack},
\begin{align*}\Bbb E(\mathbb{g}_{\u t}(\frak w_\g))&\xleftarrow[n\to\infty]{}
\int_{{\widetilde{\Om}}\cap [\v\le n]} \mathbb{g}_{\u t}(\D_{\frac{z_n}n,\g}\psi_n)dm
\\ &=
\sum_{k=1}^n\int_{\widetilde{\Om}} \mathbb{g}_{\u t}(\D_{\frac{k}n,\g}\psi_n)1_{[z_n=k]} dm\\ &=
\sum_{k=1}^n\int_{\widetilde{\Om}} \mathbb{g}_{\u t}(\tfrac{a(k)}{a(n)}(\tfrac{n}{k})^\g\psi_k)1_{{\widetilde{\Om}}\cap [\v>n-k]}\circ T^k dm\\ &\sim
\sum_{k=1}^n\int_{\widetilde{\Om}} \mathbb{g}_{\u t}(\psi_k)1_{{\widetilde{\Om}}}\circ T^k dm\cdot c(n-k)
\\ &\sim \frak e_\g(\mathbb{g}_{\u t})\cdot\sum_{k=1}^nu(k)c(n-k)\ \text{by \ref{king} on p. \pageref{king}}\
\\ &=m({\widetilde{\Om}}\cap [\v\le n])\frak e_\g(\mathbb{g}_{\u t})
\\ &\xrightarrow[n\to\infty]{×}\frak e_\g(\mathbb{g}_{\u t}) \ \ 
.\ \ \CheckedBox\ \text{ Propositions D $\&$ E.}
 \end{align*}

To establish Proposition C, recall that 
 $\psi_n\xrightarrow[n\to\infty]{\text{\tiny\tt RV}(D_{\uparrow,\infty},J_1)}\frak m_\g$,
 whence by Proposition 4 (\ref{pawn} on p. \pageref{pawn}), for $ h\in C_B(D_{\uparrow,1})$,
$$\sum_{k=1}^n\int_{\widetilde{\Om}} h(\psi_k)1_{{\widetilde{\Om}}}\circ T^k dm\sim \Bbb E(h(\text{\tt w}(\frak m_\g)))a(n)\ \text{ as $n\to\infty$};$$
whereas by Proposition E,
$$\int_{\widetilde{\Om}} h(\psi_k)1_{{\widetilde{\Om}}}\circ T^k dm\sim u(k)\Bbb E(h(\frak w_\g)),$$
and (since $a(n)\sim\sum_{k=1}^n u(k)$),
 $$\frak w_\g\ \overset{\tiny d.}=\ \text{\tt w}(\frak m_\g)$$
 proving Proposition C.\ \CheckedBox
 
 \

 We conclude this section with the:
\demo{Proof of Theorem U \ref{umbrella}}
\

By the functional Darling-Kac theorem (\cite{Bin71,Ow-Sam}), \ref{rook} (on p. \pageref{rook}) is satisfied with $\text{\tt m}_\g=\frak m_\g$. By Propositions 4 and C, \ref{bcbook} holds with $\text{\tt w}_\g(\frak m_\g)=\frak w_\g$. This proves \ref{umbrella}.\ \ \CheckedBox
\section*{\S3  Tied down occupation times}
\

In this section, we prove Theorem U \ref{knight}\ $\&$\ \ref{Industry} (p.\pageref{knight}).

\

We'll use the following   Uniform  GL Lemma, which 
generalizes \cite[Lemma 2.2.1]{G-L}.
\proclaim{Uniform  GL Lemma}
\

 For every $C=[c,d]\subset\Bbb R_+,\ \exists\ \D_n\downarrow 0$ so that  
 for $h\in L_{\a,\th},\ g\in C_B(\Bbb R),\ g, h\ge 0$, 
 on the Gibbs-Markov set $\Om$,
\begin{align*}\tag{{\tt UGL}}\label{UGL}\tfrac1{u(n)}\T^n&(hg(\tfrac{S_n(1_\Om)}{a(n)}))\ge\\ &
{m_\Om}(h)\mathbb E(g(W_\g)1_C(W_\g)))-\D_n\|g\|_{C_B}\|h\|_{L_{\a,\th}}.
\end{align*} where $W_\g$ is as in \ref{queen} on p.\pageref{queen}.
\endproclaim
We'll prove (\ref{UGL}) in the next section using ``{\tt operator renewal theoretic}'' techniques as in
\cite{OS-sub,Goue11,MelTer}.

\

We'll also use the

\

\proclaim{Gibbs Markov distortion lemma}
\

Suppose that $(Y,\mu,\tau,\a)$ is a  Gibbs-Markov map (as on p.\pageref{gmcyl}), then 
$\exists\ M>0$ so that
\begin{align*}\tag*{\PointingHand}\label{PointingHand}\widehat{\tau}^n(F)= M^{\pm 1}\|F\|_1 \ \text{\rm on}\ Y\  \forall\ n\ge 1,\ F\in L^1(Y,\s(\a_n),\mu)_+;
\end{align*}
\begin{align*}\tag*{\faMapMarker}\label{faMapMarker}\|\widehat{\tau}^n(F)\|_{L_{\a,\th}}\le M\|F\|_1 \ \ \forall\ n\ge 1,\ F\in L^1(Y,\s(\a_n),\mu). 
\end{align*}
\endproclaim
\demo{Proof}
By the Gibbs-Markov property (p.\pageref{gmcyl}), $\exists\ C>1,\ r\in (0,1)$ so that
$$|\log v_a'(x)-\log v_a'(y)|\ \le\  Cr^{t(x,y)}\ \ \ \forall\ n\ge 1,\ a\in\a_n,\ x,y\in \tau^na$$
whence 
$$v_a'(x)=C^{\pm 1}\tfrac{\mu(a)}{\mu(\tau^na)}=C^{\prime\pm 1}\mu(a)\ \forall\ a\in \a_n,\ x\in \tau^na$$
where $C'=\tfrac{C}{\min_{a\in \a}\mu(Ta)}$; and 
$$|v_a'(x)-v_a'(y)|\le CC'\mu(a)r^{t(x,y)}\ \ \ \forall\ a\in \a_n,\ x,y\in \tau^na.$$
Thus for $F:\Om\to\Bbb R\ \a_n$-measurable, non-negative,
\begin{align*}
 &\widehat{\tau}^n(F)(x)=\sum_{a\in\a_n}F(a)v_a'(x)=C^{\prime\pm 1}\int_YFd\mu\ \&\\ & \|\widehat{\tau}^n(F)\|_{L_{\a,\th}}\le \sum_{a\in\a_n}|F(a)|\|v_a'\|_{L_{\a,\th}}\le 2CC'\sum_{a\in\a_n}\mu(a)|F(a)|
\end{align*}
where $L_{\a,\th}$ is as on p.\pageref{gmcyl}. \ \CheckedBox\ 
\subsection*{Proof of Theorem U {\tt continued}}: \ \ref{knight} \ $\&$\ \ref{Industry}
\                                             

Let $\Om\in\B(X)$ be the Gibbs-Markov set with $m(\Om)=1$ and return time function $\v:\Om\to\Bbb N$. 
\

Let $(\Om,m_\Om,T_\Om,\a)$ be the corresponding Gibbs-Markov map and let
\begin{align*}
 \psi_n(t):=\tfrac{s_{[nt]}}{a(n)},\ s_n:=S_n(1_\Om),
 \ Z_n:=\min\{k\ge n+1:\ T^kx\in\Om\}.
\end{align*}
It suffices to establish \ref{knight} with $g=\mathbb{g}_{\u t}\in \Pi_0$  (as defined by \ref{faThumbTack} on p. \pageref{faThumbTack}). 
We claim that for this, it suffices to show 
\begin{align*}\tag*{\dstechnical}\label{dstechnical}
 \T^{n}(1_A\mathbb{g}_{\u t}(\psi_{n}))\ &  \gtrsim \  \Bbb E(\mathbb{g}_{\u t}(\mathfrak w_\g))m(A)u(n)\\ & 
 \text{uniformly on}\ \Om \  \forall\ A\in\mathcal C_\a\ \&\ \mathbb{g}_{\u t}\in\Pi_0
\end{align*}
where $\mathbb{g}_{\u t}\in \Pi_0$  (as defined by \ref{faThumbTack} on p. \pageref{faThumbTack}).

\

Note that \ref{dstechnical} is an extension of \ref{faFlash} on p \pageref{faFlash}
and similarly extends to hold  $\forall\ g\in\Pi_1$  (as defined by \ref{faLock} on p. \pageref{faLock}). \ 

\

Assume \ref{dstechnical} \ (whose proof is given on p.\pageref{pfdstechnical}),  then \ref{knight} with $f=1_\Om$ and $A\in\mathcal{C}_\a$ follows from \cite[Prop. 4.2]{RatWM} via
\begin{align*}\tag*{\faCog}\label{faCog}
 \tfrac1{a(N)}\sum_{n=1}^N\T^{n}(&1_A\mathbb{g}_{\u t}(\psi_{n}))\ \xrightarrow[N\to\infty]{}\  \Bbb E(\mathbb{g}_{\u t}(\mathfrak w_\g))m(A)\\ & 
 \text{a.e. on}\ \Om \  \forall\ A\in\mathcal C_\a\ \&\ \mathbb{g}_{\u t}\in\Pi_0.
\end{align*}
To extend\ \ref{knight} to hold $\forall\ f\in L^1_+\ \&\ A\in\mathcal{F}_+$, apply the argument in \cite[Proof of Theorem A:  finish]{A-Sera}.\ \ \CheckedBox\ \ref{knight}
\subsection*{Remark}
\

 Suppose that $t\in [0,1] \mapsto\mathbb{f}^{(t)}\in L^1(m)_+$ and set
$\Psi_{n,\mathbb{f}}(t):=\frac{S_{\lfl nt\rfl}(\mathbb{f}^{(t)})}{a(n)}$.
\

Variations of  the argument in \cite[Proof of Theorem A:  finish]{A-Sera} show that under the assumptions
of Theorem U \ref{umbrella} and \ref{knight} respectively, we have for $g\in\Pi$ that
\begin{align*}\tag*{{\faCloud}}
\tfrac1{a(n)}\sum_{k=1}^n&\int_Ag(\Psi_{k,\mathbb{f}})1_B\circ T^kdm\xrightarrow[n\to\infty]{}
m(A)m(B)\Bbb E(g(m(\mathbb{f})\frak w_\g))\\ & \forall\ A,\ B\in R(T)\end{align*} where $m(\mathbb{f})\frak w_\g\in\text{\tt RV}(\Bbb R_+^{[0,1]}),\ t\mapsto m(\mathbb{f}^{(t)})\frak w_\g(t)$; and
\begin{align*}\tag*{{\scriptsize\BlackKnightOnWhite}}
\tfrac1{a(n)}\sum_{k=1}^n&|\T^{k}(1_Ag(\Psi_{k,\mathbb{f}}))-u(k)m(A)\Bbb E(g(m(\mathbb{f})\mathfrak w_\g))|\\ & \xrightarrow[n\to\infty]{}\ 0\ \ \ \text{a.e.}\ 
\forall\ A\in\mathcal F_+.\end{align*}
\demo{Proof of \ref{faCog}\  given \ref{dstechnical} }
\

Let $g\in \Pi_0$, then $1-g\in \Pi_1$ and by \ref{dstechnical} (extended to $\Pi_1$) and  pointwise dual ergodicity, a.e.\ on $\Omega$,
\begin{align*}
 \mathbb{E}(g(\mathfrak w_\g))m(A)
 &\lesssim
 \tfrac1{a(N)}\sum_{n=1}^{N}\T^{n}(1_Ag(\psi_{n}))
 \\
 &=
 \tfrac1{a(N)}\sum_{n=1}^{N}\T^{n}(1_A-1_A(1-g(\psi_{n})))
 \\
 &\lesssim m(A)-\mathbb{E}(1-g(\mathfrak w_\g))m(A)
 \\
 &=\mathbb{E}(g(\mathfrak w_\g))m(A),
\end{align*}
which implies \ref{faCog} for $g\in \Pi_0$. 

\

Since $\psi_n\in D_{\uparrow,1}$ and $\mathfrak w_\g$ is continuous a.s., the convergence \ref{faCog} also holds for $g\in C_{B}(D_{\uparrow, 1})$. \ \CheckedBox

\

\demo{Proof of Theorem U\ \ref{Industry} \ given \ref{dstechnical}}

Let $g\in \Pi_0$, then $1-g\in \Pi_1$. By \ref{dstechnical}\ and (\ref{OSRT}), for $A\in\mathcal{C}_\a$, uniformly \ on $\Omega$,
\begin{align*}
 \mathbb{E}(g(\mathfrak w_\g))m(A)
 &\lesssim
 \tfrac1{u(n)}\T^{n}(1_Ag(\psi_{n}))
 \\
 &=
 \tfrac1{u(n)}\T^{n}(1_A)-\tfrac1{u(n)}\T^{n}(1_A(1-g(\psi_{n})))
 \\
 &\lesssim m(A)-\mathbb{E}(1-g(\mathfrak w_\g))m(A)
 \\
 &=\mathbb{E}(g(\mathfrak w_\g))m(A).
\end{align*}
Therefore we obtain the desired uniform convergence for $g\in \Pi_0$. This convergence also holds for $g\in C_B(D_{\uparrow, 1})$ since $\mathfrak w_\g$ is a continuous process. 
\CheckedBox

\demo{\bf Proof of \ref{dstechnical}}\label{pfdstechnical}\footnote{ on p. \pageref{dstechnical}}

The rest of this section is devoted to a proof of  \ref{dstechnical} 
given the uniform GL lemma which enables a suitable
 modification of the proof of \ref{faFlash} (p. \pageref{faFlash}) 
 without using a Markov property.

\

Again writing $\kappa:=\lfl nt_N\rfl$, fix $A\in \a_J$, then
\begin{align*}1_A \mathbb{g}_{\u t}(\psi_n)&=\sum_{k=\kappa+1}^n1_A \mathbb{g}_{\u t}(\psi_n)1_{[Z_\kappa=k]}\\ &=
\sum_{k=\kappa+1}^n1_A\frak g(\psi_n)1_{[Z_\kappa=k]}1_\Om\circ T^k h(\psi_{n}(1)-\psi_n(t_N))\\ &=
\sum_{k=\kappa+1}^n1_A\frak g(\psi_n)1_{[Z_\kappa=k]}1_\Om\circ T^k h(\tfrac{s_{n-k}\circ T^k}{a(n)});
\end{align*}
and 
                                              
\begin{align*}\tag*{\dsheraldical}\label{dsheraldical}\T^{n}(1_A&\mathbb{g}_{\u t}(\psi_n))=
\\ &=\sum_{k=\kappa+1}^n\T^n\big(1_A\frak g(\psi_n)1_{[Z_\kappa=k]}1_\Om\circ T^k h(\tfrac{a(n-k)}{a(n)}\tfrac{s_{n-k}\circ T^k}{a(n-k)})\big)
\\ &=\sum_{k=\kappa+1}^n\T^{n-k}\big(\T^k(1_A\frak g(\psi_n)1_{[Z_\kappa=k]})1_\Om h(\tfrac{a(n-k)}{a(n)}\tfrac{s_{n-k}}{a(n-k)})\big).\end{align*}

Next, for fixed $ \kappa+1\le k\le n$, define on $\Om$,
\begin{align*}G_{k,n}&:=\T^k(1_A\frak g(\psi_n)1_{[Z_\kappa=k]})\\ &=
\sum_{\ell=1}^k\T_\Om^\ell(1_A\frak g(\psi_n)1_{[\varphi_\ell=Z_\kappa=k]})
\end{align*} where $\v_\ell:=\sum_{k=0}^{\ell-1}\v\circ T_\Om^j$.
\

We claim that $\exists\ C^{\pprime}>0$ so that
\begin{align*}\tag*{{\scriptsize\faPaw}}\label{faPaw}\|G_{k,n}\|_{L_{\a,\th}}\le C^{\pprime}\|\frak g\|_{C_B} m_\Om([Z_\kappa=k])\ \forall\  \kappa+1\le k\le n.
\end{align*}
\demo{Proof of \ref{faPaw}}\ \ 
For $\kappa$ as above and $\kappa+1\le k\le n,\ \ell\le k$, let 
$$F_{k,\ell,n}:=1_A\frak g(\psi_n)1_{[\varphi_\ell=Z_\kappa=k]}$$ 

is $\a_\ell$-measurable and by \ref{faMapMarker} (p. \pageref{faMapMarker}), $G_{k,\ell,n}:=\T_\Om^\ell F_{k,\ell,n}\in L_{\a,\th}$ (as on p.\pageref{gmcyl}) with
$$\|G_{k,\ell,n}\|_{L_{\a,\th}}\le C^{\pprime} m_\Om(|F_{k,\ell,n}|)\le 
C^{\pprime} \|\frak g\|_{C_B} m_\Om([\v_\ell=Z_\kappa=k]).$$
where $\|\frak g\|_{C_B}:=\sup_{\xi\in D_{\uparrow,R}}|\frak g(\xi)|$.
Thus
\begin{align*}
\|G_{k,n}\|_{L_{\a,\th}}&\le\sum_{\ell=1}^k\|G_{k,\ell,n}\|_{L_{\a,\th}}\le 
C^{\pprime}\sum_{\ell=1}^k \|\frak g\|_{C_B} m_\Om([\v_\ell=Z_\kappa=k])\\ &=C^{\pprime}\|\frak g\|_{C_B} m_\Om([Z_\kappa=k]).\ \ \CheckedBox\ \ \text{\ref{faPaw}}
\end{align*}

Fix $C=[c,d]\subset\Bbb R_+$.  By (\ref{UGL}) and \ref{faPaw},
\begin{align*}\tfrac1{u({n-k})}&\T^{n-k}(G_{n,k}1_\Om h(\tfrac{a(n-k)}{a(n)}\tfrac{s_{{n-k}}}{a(n-k)})\big)\\ &\ge m_\Om(G_{n,k})\Bbb E(h(\tfrac{a(n-k)}{a(n)}W_\g) 1_C(W_\g))-\D_{n-k}  \|G_{n,k}\|_{L_{\a,\th}}\|h\|_{C_B}\\ &=
\int_A\frak g(\psi_n)1_{[z_\kappa=k]}dm\Bbb E(h(\tfrac{a(n-k)}{a(n)})W_\g)1_C(W_\g))-\D_{n-k}    C^{\pprime} \|g\|_{C_B} m_\Om([Z_\kappa=k])
.\end{align*}

It follows  that for fixed $t_N<\l<1$, 
\begin{align*}\T^{n}(1_A&\mathbb{g}_{\u t}(\psi_n))\ge\\ &
\sum_{\kappa+1\le k\le \l n}\int_A\frak g(\psi_n)1_{[Z_\kappa=k]}dm\cdot \Bbb E(h(\tfrac{a(n-k)}{a(n)} W_\g)1_C(W_\g))u(n-k)\\ &\ \ \ \ \ 
-  C^{\pprime} \|g\|_{C_B}\sum_{\kappa+1\le k\le \l n}m_\Om([Z_\kappa=k])u(n-k)\D_{n-k}.\end{align*}
To continue, fix $\eta>0$ and let $J=J(C,\l,\eta)$ be so that $\forall\ n\ge J,\  \kappa+1\le k\le \l n,$
\begin{align*}u(n-k)=(1\pm\eta)\tfrac{u(n)}{(1-\frac{k}n)^{1-\g}}\ \&\ h(t\tfrac{a(n-k)}{a(n)})=(1\pm\eta)h(t(\tfrac{n-k}n)^\g)\ \forall\ t\in C. 
\end{align*}

For $n\ge J$,
\begin{align*}\sum_{\kappa+1\le k\le \l n}m_\Om([Z_\kappa=k])u(n-k)\D_{n-k}&\le (1+\eta)u(n)\sum_{\kappa+1\le k\le \l n}m_\Om([Z_\kappa=k])(\tfrac{n}{n-k})^{1-\g}\D_{n-k}\\ &
\le (1+\eta)u(n)(\tfrac1{1-\l})^{1-\g}\D_{(1-\l)n} 
\end{align*}
and
\begin{align*}
\sum_{\kappa+1\le k\le \l n}&\int_A\frak g(\psi_n)1_{[Z_\kappa=k]}dm\cdot \Bbb E(h(\tfrac{a(n-k)}{a(n)} W_\g) 1_C(W_\g))u(n-k)\\ &\ge
(1-\eta)^2u(n)\sum_{\kappa+1\le k\le \l n}\int_A\frak g(\psi_n)1_{[Z_\kappa=k]}dm\cdot \Bbb E(h((\tfrac{n-k}n)^\g W_\g) 1_C(W_\g))(\tfrac{n}{n-k})^{1-\g}.
\end{align*}
Moreover, with $(\psi_n,Z_\kappa)\perp W_\g$,
\begin{align*}
\sum_{\kappa+1\le k\le \l n}&\int_A\frak g(\psi_n)1_{[Z_\kappa=k]}dm\cdot \Bbb E(h((\tfrac{n-k}n)^\g  W_\g)1_C(W_\g))(\tfrac{n}{n-k})^{1-\g}\\ &=\int_A\frak g(\psi_n)\sum_{\kappa+1\le k\le \l n}\(1_{[Z_\kappa=k]}\cdot \Bbb E(h((\tfrac{n-k}n)^\g  W_\g)1_C(W_\g))(\tfrac{n}{n-k})^{1-\g}\)dm\\ &=
\int_A\(\frak g(\psi_n)1_{[Z_\kappa\le\l n]}\Bbb E(h(\tfrac{{n-Z_\kappa}}n)^\g W_\g1_C(W_\g)\|Z_\kappa)(\tfrac{n}{n-Z_\kappa})^{1-\g}\)dm\\ &=
\int_A\(\frak g(\psi_n)1_{[Z_\kappa\le\l n]}\Bbb E(h((\tfrac{{n-Z_\kappa}}n)^\g  W_\g)1_C(W_\g))(\tfrac{n}{n-Z_\kappa})^{1-\g}\)dm \ \because\ (\psi_n,Z_\kappa)\perp W_\g.
\end{align*}

\

Using \ref{dscommercial} on p.\pageref{dscommercial} and   Eagleson's theorem (\cite{Eagleson}),
{\small\begin{align*}
\int_A&\(\frak g(\psi_n)1_{[Z_\kappa\le\l n]}\Bbb E(h((\tfrac{{n-Z_\kappa}}n)^\g  W_\g)1_C(W_\g))(\tfrac{n}{n-Z_\kappa})^{1-\g}\)dm
\xrightarrow[n\to\infty]{}m(A)\cdot\\ &\cdot\Bbb E\(\frak g(\frak m_\g)1_{[\mathcal D(\frak m_\g)(t_N)\le \l]}(\tfrac1{1-\mathcal D(\frak m_\g)(t_N)})^{1-\g}
h((1-\mathcal D(\frak m_\g)(t_N))^\g  W_\g)1_C(W_\g)\)\\ &
\xrightarrow[\l\to 1-,\ C\uparrow\Bbb R_+]{}\ m(A)\frak e_\g(\mathbb{g}_{\u t})\\ &= m(A)\Bbb E(\mathbb{g}_{\u t}(\frak w_\g))\ \ \text{by Proposition D on p. \pageref{queen}}.
\end{align*}}
Thus,
\begin{align*}\T^{n}(1_A\mathbb{g}_{\u t}(\psi_n))
\gtrsim 
m(A)\Bbb E(\mathbb{g}_{\u t}(\frak w_\g)) u(n)\ \ \text{uniformly on $\Om$.\ \ \Checkedbox\ \ \ref{dstechnical}}.
\end{align*}

\section*{\S4  Uniform {\tt LLT}s $\&$ the uniform GL lemma}

In this section, we'll prove the uniform GL lemma \ref{UGL} (p. \pageref{UGL}).
Throughout the section, we use notations from \cite[\S3]{A-Sera} 

\

Let $(\Om,\mu,\tau,\a)$ be a    mixing, probability
preserving Gibbs-Markov map as on p.\pageref{gmcyl} and let
$\phi :\Om\to \Bbb N$
be 
locally $L_{\a,\th}$,  in the sense that  $\exists\ M>0$ so that $$\forall\ A\in\a,\ 1_A\phi\in L_{\a,\th}\  \&\ 
\|1_A\phi-\int_A\phi d\mu\|_{L_{\a,\th}}\le M;$$
\ and also satisfy 
\begin{align*}\tag*{{\scriptsize\faKey}}\label{faKey}
 \mu([\phi>t])\  \sim\ \frac1{\G(1+\g)\G(1-\g)}\cdot\frac1{a(t)}.
\end{align*}
 with $a(t)\ \g$-regularly varying with $0<\g<1$.

 \

 As in Nagaev's theorem  (\cite{Nagaev}, see also  \cite[Theorem 4.1]{AD} $\&$ the proof of Theorem 3.1 in \cite{A-Sera}), we write:
$$P_{\phi,t}(f):=\widehat{\tau}(e(t\phi)f),\ N(t)=N_\phi(t):=N(P_{\phi,t})\ \ ,\ \ \l(t)=\l_\phi(t):=\l(P_{\phi,t}).$$
(with $e(t):=e^{2\pi it}$).
\

By \cite[theorem 6.1]{AD}, under the above conditions,
$$\tfrac{\phi_n}{a^{-1}(n)}\xrightarrow[n\to\infty]{\frak d} \ Z_\g:= \eta_\g(1)$$
where $\phi_n:=\sum_{k=0}^{n-1}\phi\circ\tau^k$;  and
$$\l(\tfrac{t}{a^{-1}(n)})^n\xrightarrow[n\to\infty]{}\widehat{f}_{Z_\g}(t).$$

\

We'll need a uniform  version of the  lattice {\tt LLT}.

\
As in \cite[\S3]{A-Sera}, suppose that $\phi:\Om\to\Bbb N$ is non-arithmetic, then either $\phi$ is aperiodic, or
$\phi$ has a reduction
\begin{align*}\tag*{\phone}\label{phone}
 \phi=\ F-F\circ \tau+{p}\psi+\xi
\end{align*}

where  $F:\Om\to\Bbb N,\ F\in L_{\a,\th}$, ${p}\in\Bbb N,\ p\ge 2,\ \xi\in [1,{p}-1]\cap\Bbb N,\ \text{\tt gcd}\{\xi,{p}\}=1$  and $\psi:\Om\to \Bbb Z$  is locally $L_{\a,\th}$  and aperiodic.
\

\proclaim{Uniform  conditional lattice {\tt LLT}}
\par Under the above conditions, 
$\exists\ \D(N)\downarrow 0$ so that
if 
$$k_n\in\Bbb Z,\ \kappa_n:=\tfrac{k_n}{a^{-1}(n)}\ \forall\ n\ge 1,$$
then $\forall\  n\ge 1,\ h\in L_{\a,\th}$:
\

in case $\phi$ is aperiodic:
 \begin{align*}\tag*{\Tumbler}\label{Tumbler}
  \|a^{-1}(n)\ttau^n(h1_{[\phi_n=k_n]})-\mu(h)f_{Z_\g }(\kappa_n)\|_{L_{\a,\th}}\le \D(n)\|h\|_{L_{\a,\th}};
 \end{align*}
and, in case $\phi$ has the reduction \ref{phone},
 {\small\begin{align*}\tag*{\Bleech}\label{Bleech}\|a^{-1}(n)\widehat{\tau}^n(h1_{[\phi_n=k_n]})-1_{p\Bbb Z}(n\xi-k_n)f_{Z_\g}(\kappa_n)\mu(h)\|_{L_{\a,\th}}\le \D(n)\|h\|_{L_{\a,\th}}.
 \end{align*}}
\endproclaim
\f{\bf Remarks}
\

\ \ The estimates are useful when $(\kappa_n:\ n\ge 1)$ is bounded away from $0\ \&\ \infty$ (i.e.  $\inf_nf_{Z_\g}(\kappa_n)>0$). 
\

   For  estimates outside this range,  see the ``local large deviations'' of \cite{C-D,MelTer20} and the ``extended  {\tt LLT}s'' of \cite{A-T} and references therein.
\ 

The periodic {\tt LLT} \ref{Bleech} was first proved in \cite{Shepp}
(in the independent case).
\
 
\demo{(i) Proof of \ref{Tumbler}}
By \cite[Thm. 4.1]{AD}, 
$\exists\ \d>0,\ \th\in (0,1),\ C>1$ such that $\forall\  n\ge 1,\ h\in {L_{\a,\th}}$:
$$\|P_{\phi,t}^nh-\l(t)^nN(t)h\|_{L_{\a,\th}}\le C\|h\|_{L_{\a,\th}}\th^n\ \ \forall\ |t|\le\d,$$
 and 
$$\|P_{\phi,y}^n1\|_{L_{\a,\th}}\le C\|h\|_{L_{\a,\th}}\th^n\ \ \forall\ \d\le |y|\le \pi.$$
Write $a(t)=\tfrac{t^\g}{\ell(t)}$ with $\ell(t)$ slowly varying at infinity.
\

Using  \cite[Thm. 5.1]{AD}, by possibly
shrinking $\d>0$, we can ensure in addition that
\begin{align*}\tag*{\faBomb}\label{faBomb}-\text{Re}\,\log\l(t)
\ge \frac{c}{2a(\frac1t)}\ \ \forall\ |t|\le\d
 \end{align*}

and by slow variation of $\ell$,
$\exists\ 0<\e=\e(\d)$ such that
$${\ell({a^{-1}(n)\over |t|})\over \ell(a^{-1}(n))}\ge |t|^\e\ \forall\ n\ge 1,\ |t|\le\d a^{-1}(n).$$
   The integrals below converge in  $L_{\a,\th}$ and we have,  in $L_{\a,\th}$,
\begin{align*} \ 2\pi a^{-1}(n)\ttau^n(h1_{[\phi_n=k_n]}) &=
a^{-1}(n)\ttau^n\(h\int_{-\pi}^{\pi}e^{-itk_n}e^{it\phi_n}dt\)\\ &=
a^{-1}(n)\int_{-\pi}^{\pi}e^{-itk_n}\ttau^n(he^{it\phi_n})dt
\\ &=a^{-1}(n)\int_{-\pi}^{\pi}e^{-itk_n}P_{\phi,t}^nhdt\\ &=a^{-1}(n)\int_{|t|\le\d}e^{-itk_n}\l(t)^nN(t)hdt\pm 2\pi C\|h\|_{L_{\a,\th}}a^{-1}(n)\th^n.\end{align*}
Next
\begin{align*}
a^{-1}(n)\int_{|t|\le\d}&e^{-itk_n}\l(t)^nN(t)hdt =
\int_{-\d a^{-1}(n)}^{\d a^{-1}(n)}e^{-it\kappa_n}\l(\tfrac{t}{a^{-1}(n)})^n
N(\tfrac{t}{a^{-1}(n)})hdt \\ &\&\ \ \ \ \ \ \ 2\pi f_{Z_\g }(\kappa)\ =\ 
\int_{\Bbb R}\widehat{f}_{Z_\g }(t)e^{-i\kappa t}dt,\end{align*}
so
\begin{align*}
a^{-1}(n)\int_{|t|\le\d}&e^{-itk_n}\l(t)^nN(t)hdt -2\pi \mu(h) f_{Z_\g }(\kappa_n)\\ &=
\int_{|t|\le \d a^{-1}(n)}e^{-it\kappa_n}(\l(\tfrac{t}{a^{-1}(n)})^n
N(\tfrac{t}{a^{-1}(n)})h-\mu(h)\widehat{f}_{Z_\g }(t))dt\\ &\ \ \ \ \ \ \ \ \ \ +\int_{|t|>\d a^{-1}(n)}e^{-it\kappa_n} \widehat{f}_{Z_\g }(t)dt\mu(h)\\ &=:I_n+II_n.\end{align*}
Now $|\widehat{f}_{Z_\g}(t)|=e^{-c|t|^\g}$ ($c>0$) so
    $$\|II_n\|_{L_{\a,\th}}\le \|h\|_{L_{\a,\th}}\int_{|t|>\d a^{-1}(n)} e^{-c|t|^\g}dt\xrightarrow[n\to\infty]{}\ 0.$$
By \cite[Theorems 2.4 $\&$ 5.1]{AD},
\begin{align*}\tag*{{\scriptsize\faTaxi}}
 \|N(\tfrac{t}{a^{-1}(n)})h-\mu(h)\|_{L_{\a,\th}}&\le Mm_\Om(|e^{i\tfrac{t\phi}{a^{-1}(n)}}-1|)
 \le M''\tfrac1{a(\frac{a^{-1}(n)}t)}\|h\|_{L_{\a,\th}}\\ &\sim M''\tfrac{t^{\frac1\g}}{a(a^{-1}(n))}\|h\|_{L_{\a,\th}}\sim M''\tfrac{t^{\frac1\g}\|h\|_{L_{\a,\th}}}n
\end{align*}
and 
\begin{align*}\|I_n\|_{L_{\a,\th}} &=\int_{-\d a^{-1}(n)}^{\d a^{-1}(n)}\|\l(\tfrac{t}{a^{-1}(n)})^n
N(\tfrac{t}{a^{-1}(n)})h-\mu(h)\widehat{f}_{Z_\g }(t)\|_{L_{\a,\th}}dt\\ &\le
\int_{-\d a^{-1}(n)}^{\d a^{-1}(n)}|\l(\tfrac{t}{a^{-1}(n)})|^n\|N(\tfrac{t}{a^{-1}(n)})h-\mu(h)\|_{L_{\a,\th}}dt\\ &\ \ \ \ \ \ \ \ \ \ \ \ \ \ \ \ \ +|\mu(h)\int_{-\d a^{-1}(n)}^{\d a^{-1}(n)}|\l(\tfrac{t}{a^{-1}(n)})^n
-\widehat{f}_{Z_\g }(t)|dt\\ &\overset{\text{\tiny(\faTaxi)}}\le \|h\|_{L_{\a,\th}}\left(\tfrac{M''}n\int_{-\d a^{-1}(n)}^{\d a^{-1}(n)}|t|^{\frac1\g}|\l(\tfrac{t}{a^{-1}(n)})|^ndt+
\int_{-\d a^{-1}(n)}^{\d a^{-1}(n)}|\l(\tfrac{t}{a^{-1}(n)})^n
-\widehat{f}_{Z_\g }(t)|dt\right).
\end{align*}
By \ref{faBomb} (on p. \pageref{faBomb}),and regular variation of $a^{-1},\ \exists\ \e>0$ so that $|\l(\tfrac{t}{a^{-1}(n)})|^n\ll   e^{-\e|t|^{\frac{\g}2}}$. 
\

Moreover,  
$|\widehat{f}_{Z_\g }(t)|\ll e^{-c|t|^\g}$ so 
\begin{align*}&\int_{-\d a^{-1}(n)}^{\d a^{-1}(n)}|t|^{\frac1\g}|\l(\tfrac{t}{a^{-1}(n)})|^ndt\le M\int_\Bbb R|t|^{\frac1\g}e^{-\e|t|^{\frac{\g}2}}dt<\infty,\ \text{whence}
\\ & \int_{-\d a^{-1}(n)}^{\d a^{-1}(n)}|\l(\tfrac{t}{a^{-1}(n)})^n
-\widehat{f}_{Z_\g }(t)|dt\xrightarrow[n\to\infty]{}\ 0\\ &\ \ \  \text{by dominated convergence.\ \ \Checkedbox\ \ \ref{Tumbler}}.
\end{align*}

\demo{(ii) Proof of \ref{Bleech}}\ \footnote{This proof is an modification of the proof of  the lattice case of \cite[Thm. 3.1]{A-Sera}.}
\

Recall from \ref{phone} on p. \pageref{phone} that
$\phi=p\psi+\xi+F-F\circ\tau$ where
\sbul $p,\ \xi\in\Bbb N,\ 1\le\xi<p,\ \text{\tt gcd}\,\{p,\xi\}=1$;
\sbul  $\psi:\Om\to\Bbb N$ is aperiodic $\&\ F:\Om\to\Bbb Z$ satisfies $F\circ\tau-F\in p\Bbb Z$.
\

For $k_n\in\Bbb Z$,
$$\phi_n=k_n\ \iff\ p\psi_n=k_n-n\xi-F+F\circ\tau^n$$
and this entails $k_n-n\xi\in p\Bbb Z$.
\

Thus,
 \

  \begin{align*}\widehat{\tau}^n(&h1_{[\phi_n=k_n]})(x)=1_{{p}\Bbb Z}(n\xi-k_n)\int_\Bbb T e(-k_nt)\ttau^n(e(t\phi_n)h)(x)dt\\ &=
    1_{{p}\Bbb Z}(n\xi-k_n)\int_\Bbb T e(t(n\xi-k_n))\ttau^n(e(pt\psi_n)e(t(F-F\circ\tau^n))h)(x)dt\\ & =
    1_{{p}\Bbb Z}(n\xi-k_n)\int_\Bbb T e(pt\cdot\tfrac{n\xi-k_n}p)\ttau^n(e(pt\psi_n)e(pt\cdot\tfrac{F-F\circ\tau^n}p))h)(x)dt\\ & =
        {p}1_{{p}\Bbb Z}(n\xi-k_n)\int_\Bbb T e(\tfrac{t}{p}(n\xi-k_n))\ttau^n(e(t\psi_n)e(\tfrac{t}{p}(F-F\circ\tau^n))h)(x)dt\\ &\ \ \ \ \ \ \ \ \ \ \ \  \ \because\ t\mapsto\ pt\ \text{is measure preserving on}\ \Bbb T;\\ &=
      {p}1_{{p}\Bbb Z}(n\xi-k_n)\int_\Bbb T e(\tfrac{t}{p}(n\xi-k_n-F(x)))P_{\psi,t}^n(e(\tfrac{tF}{p})h)(x)dt.
   \end{align*}
   By the proof of \ref{Tumbler},
\begin{align*}\|a^{-1}(n)\int_\Bbb T  e(\tfrac{t}{p}(n\xi-k_n-F(x)))P_{\psi,t}^n(e(\tfrac{tF}{p})h)(x)dt-\mu(h)f_{\frac{Z_\g }{p}}(\tfrac{\kappa_n}{p})\|_{L_{\a,\th}} \le \D(n)\|h\|_{L_{\a,\th}}
\end{align*}
and \ref{Bleech} follows from this.\ \Checkedbox

\demo{Proof of the uniform GL lemma}

As in the proof of \cite[Lemma 2.1(GL)]{A-Sera}, we assume WLOG 
 that $a^{-1}(n+1)-a^{-1}(n)\sim\tfrac{a^{-1}(n)}{\g n}$ and set $x_{\nu ,n}:=\tfrac{n}{a^{-1}(\nu )}$; obtaining (see \cite{A-Sera})
 \begin{align*}\tag*{{\scriptsize\faLinux}}
 \tfrac1{a^{-1}(k)}\sim\tfrac{\g k}n\cdot (x_{k,n}-x_{k+1,n})\sim \tfrac{\g a(n)}n\cdot\tfrac{x_{k,n}-x_{k+1,n}}{x_{k,n}^\g}=u(n)\cdot\tfrac{x_{k,n}-x_{k+1,n}}{x_{k,n}^\g}.
\end{align*}

It suffices to consider $g\in C_B(\Bbb R_+)$ so that $\log G:[c,d]\to\Bbb R$ is smooth where
$G(x):=g(\tfrac1{x^\g})f(x)$. 
\

Fixing $h\in L_{T_\Om},\ h\ge 0$,  we have that
\begin{align*}
&\frac1{u(n)}\T^n(hg(\tfrac{S_n(1_\Om)}{a(n)}))=\sum_{\nu =1}^n\T_\Om^{\nu }(h1_{[\v_\nu =n]}g(\tfrac{\nu }{a(n)}))\\ &\ge
\frac1{u(n)}\sum_{1\le \nu \le n,\ x_{\nu ,n}\in [c,d]}\T_\Om^{\nu }(h1_{[\v_\nu =n]}g(\tfrac{\nu }{a(n)}))\\ &\ge
\frac1{u(n)}\sum_{1\le \nu \le n,\ x_{\nu ,n}\in [c,d]}\tfrac1{a^{-1}(\nu )}[1_{p\Bbb Z}(n -\nu \xi)m_\Om(h)g(\tfrac{\nu }{a(n)})f_{Z_\g}(x_{\nu,n} )-\|g\|_{C_B}\|h\|_{L_{T_\Om}}\D(\nu) ]\ \text{by \ref{Bleech}}\\ &=
m_\Om(h)\frac1{u(n)}\sum_{1\le \nu \le n,\ x_{\nu ,n}\in [c,d]}\tfrac1{a^{-1}(\nu )}1_{p\Bbb Z}(n -\nu \xi)g(\tfrac{\nu }{a(n)})f_{Z_\g}(x_{\nu,n})\ -\|g\|_{C_B}\|h\|_{L_{T_\Om}}\mathcal E^{(1)}_n\end{align*}
where
\begin{align*}\mathcal E^{(1)}_n:=\frac1{u(n)}\sum_{1\le \nu \le n,\ x_{\nu ,n}\in [c,d]}\tfrac{\D(\nu)}{a^{-1}(\nu )}
\end{align*}
with $\D(\nu)$ as in \ref{Tumbler} on p. \pageref{Tumbler}.
\

To see that $\mathcal E^{(1)}_n\xrightarrow[n\to\infty]{}0$, we note that for  $1\le \nu \le n$,
$$x_{\nu ,n}\in [c,d]\ \ \Leftrightarrow\ \nu\in [a(\tfrac{n}d),a(\tfrac{n}c)].$$
Thus 
\begin{align*}\mathcal E^{(1)}_n&=\tfrac1{u(n)}\sum_{1\le \nu \le n,\ x_{\nu ,n}\in [c,d]}\tfrac{\D(\nu)}{a^{-1}(\nu )}
\\ &\le \tfrac{\D(a(\frac{n}d))}{u(n)}\sum_{\nu\ge a(\tfrac{n}d)}\tfrac1{a^{-1}(\nu )}\\ &
\asymp \tfrac{\D(a(\frac{n}d))}{u(n)}\tfrac{a(\frac{n}d)}{a^{-1}(a(\frac{n}d))}\\ &\asymp \D(a(\tfrac{n}d))\xrightarrow[n\to\infty]{}0. 
\end{align*}
\footnote{Here, for $a_n, b_n>0,\ a_n\asymp b_n$ means $a_n\ll b_n\ \&\ b_n\ll a_n$ where $\ll$ is as in the footnote on p. \pageref{ll}.}

\

Next, again as in the proof of \cite[Lemma 2.1(GL)]{A-Sera}, 
\begin{align*}\frac1{u(n)}&\sum_{1\le \nu \le n,\ x_{\nu ,n}\in [c,d]}\tfrac1{a^{-1}(\nu )}1_{p\Bbb Z}(n -\nu \xi)g(\tfrac{\nu }{a(n)})f_{Z_\g}(x_{\nu,n}) \\ &=
\frac1{u(n)}\sum_{1\le k\le n,\ x_{k,n}\in [c,d]}g(\tfrac1{x_{k,n}^\g})\tfrac{pf_{Z_\g}(x_{k,n})}{a^{-1}(k)}1_{p\Bbb Z}(n-k\xi)
 \\ &=  \Bbb E(1_{[c,d]}(Z_\g)g(Z_\g^{-\g})Z_\g^{-\g})\ \pm\ \mathcal E^{(2)}_n\|g\|_{C_B}
\end{align*}
where $\mathcal E^{(2)}_n\to 0$.

Thus,
\begin{align*}\tfrac1{u(n)}&\T^n(hg(\tfrac{S_n(1_\Om)}{a(n)}))\ge m_\Om(h)\mathbb E(g(W_\g)1_C(W_\g)))-\\ &-(\mathcal E^{(1)}_n+\mathcal E^{(2)}_n)\|g\|_{C_B}\|h\|_{L_{\a,\th}}.\ \CheckedBox
\end{align*}

\section*{\S5 Remarks on  $1$-self similar limit processes with  stationary increments}
\

Let $(\Om,\mu,\tau)$ be a probability preserving
transformation  $\&$ let $\v:X\to\Bbb R_+$ be measurable and suppose that $\tfrac{\v_n}{b(n)}\xrightarrow[n\to\infty]{\text{\tt\tiny RV}(\Bbb R_+)}Y\in\text{\tt RV}(\Bbb R_+)$ 
with $b:\Bbb R_+\to\Bbb R_+$ is $1$-regularly varying.
\

As remarked in \cite{AWdistlim},
\begin{align*}\tag*{{\scriptsize\faHome}}\label{faHome}\Phi_n\xrightarrow[n\to\infty]{\frak d}Y\mathbb{i}\end{align*}
in $(D_{\uparrow,\infty},\text{\tt DF})$
where 
$\Phi_n(t):=\frac{\v_{\lfl nt\rfl}}{b(n)}\ \&\ \mathbb{i}(t):=t$.
\

This follows from the $1$-regular variation of $b:\Bbb R_+\to\Bbb R_+$ via  \cite[Theorem 3.1]{vervaat85}.

Now let $(X,m,T)$ be a conservative, ergodic,  measure preserving 
transformation  and suppose that 
\begin{align*}\tfrac{S_n(f)}{a(n)}\xrightarrow[n\to\infty]{\frak d}Zm(f)\ \forall\ f\in L^1(m)_+ 
\end{align*}
where $Z\in\text{\tt RV}(\Bbb R_+)$ and  $a:\Bbb R_+\to\Bbb R_+$ is $1$-regularly varying.
We claim that
\begin{align*}\tag*{{\scriptsize\faUniversity}}\label{faUniversity}\Psi_{n,f}\xrightarrow[n\to\infty]{\frak d}m(f)Z\mathbb{i}\end{align*}
in $(D_{\uparrow,\infty},\text{\tt DF})$ where 
$\Psi_{n,f}(t):=\frac{S_{\lfl nt\rfl}}{a(n)}$.
\endproclaim\demo{Proof of ({\scriptsize\faUniversity})}\ \ Fix $\Om\in\B(X),\ m(\Om)=1$ and let
$\mu:=m_\Om,\ \v:\Om\to\Bbb N,\ \v(x):=\min\,\{n\ge 1:\ T^nx\in \Om\}\ \&\ \tau:\Om\to\Om,\ \tau(x):=T^{\v(x)}(x)$.
\

By inversion $\tfrac{\v_n}{b(n)}\xrightarrow[n\to\infty]{\frak d}Z^{-1}\in\text{\tt RV}(\Bbb R_+)$ with $b:=a^{-1}$
and by ({\scriptsize\faHome}),
\begin{align*}\Phi_n\xrightarrow[n\to\infty]{\frak d}Z^{-1}\mathbb{i}\ \text{in $(D_{\uparrow,\infty},\text{\tt DF})$},\end{align*}
whence, again by inversion
$$\Psi_{n,1_\Om}\xrightarrow[n\to\infty]{\frak d}Z\mathbb{i}\ \text{in $(D_{\uparrow,\infty},\text{\tt DF})$}.$$

Conclude to obtain ({\scriptsize\faUniversity}) by applying the ratio ergodic theorem and the functional version of Eagleson's theorem (\cite{Eagleson,TZ}.\ \Checkedbox
\

Examples of conservative, ergodic,  measure preserving 
transformations satisfying ({\scriptsize\faUniversity}) are given in \cite{AWdistlim}.

\proclaim{Proposition 6}\ \ Suppose that $(X,m,T)$ is a weakly rationally ergodic measure preserving
transformation \ \footnote{as in \ref{faTree} on p.\pageref{faTree}}
so that
\begin{align*}\tag*{{\scriptsize\faCab}}\label{faCab} \tfrac{S_n(f)}{a(n)}\xrightarrow[n\to\infty]{m} 
 m(f)\ \forall\ f\in L^1(m),\end{align*} then
\begin{align*}\tag*{\dschemical}\label{dschemical} \tfrac1{a(n)}&
\sum_{k=1}^n|\int_Ag(\tfrac{S_k(f)}{a(k)})1_B\circ T^kdm-g(m(f))m(A\cap T^{-k}B)|\\ &
\xrightarrow[n\to\infty]{}0\ \forall\ f\in L^1(m)_+,\ g\in C_B(\Bbb R_+)\ ,\ A,B\in R(T).\end{align*}
If, in addition, $(X,m,T)$ is  {\tt RWM}, then
\begin{align*}\tag*{\scriptsize\faTrain}\label{faTrain}\tfrac1{a(N)}&\sum_{n=1}^N|\int_{B\cap T^{-n}C}g(\tfrac{S_n(f)}{a(n)})dm-m(B)m(C)
 g(m(f))u_n|\\ &
\xrightarrow[N\to\infty]{}\ 0\ \forall\ B,C\in  R(T),\ g\in C_B(\Bbb R_+)\ \&\ f\in L^1_+
\end{align*} where $a(n)\sim\sum_{k=1}^nu_k.$
\endproclaim
Note that   \ref{dschemical} does not entail spectral weak mixing; and 
\ref{faTrain} is  the integrated, tied-down renewal mixing property in \cite[Theorem 6.2]{A-Sera} with $\g=1\ \&\ W_1\equiv 1$ when $u_n\sim\tfrac{a(n)}n$.
\demo{Proof of\ref{dschemical}}\ \ By \cite[Theorem 4.1]{homog}, $a(n)$ is $1$-regularly varying. Thus by Corollary 3 and Proposition 4,
\begin{align*} 
\tfrac1{a(n)}\sum_{k=1}^n&\int_Ag(\tfrac{S_k(f)}{a(k)})1_B\circ T^kdm\xrightarrow[n\to\infty]{}
m(A)m(B)g(m(f))\\ & \forall\ A,\ B\in R(T),\ f\in L^1_+,\ g\ \text{\tt Riemann integrable on}\ \Bbb R_+.\end{align*}
Fix $f\in L^1(m)_+,\ \e\in (0,1)$ and let $g_\e:=1_{((1-\e)m(f),(1+\e)m(f))}$, then writing $s_n:=S_n(f)$, we have
$$\tfrac1{a(n)}\sum_{k=1}^n\int_Ag_\e(\tfrac{s_k}{a(k)})1_B\circ T^kdm\xrightarrow[n\to\infty]{}
m(A)m(B)\ \forall\ A,\ B\in R(T).$$
For $\e>0\ \&\ x\in X$, let
$$K(\e,x):=\{k\ge 1:\ |\tfrac{s_k}{a(k)}-m(f)|<\e m(f)\}.$$ By the above, for $A,\ B\in R(T)$,
$$\int_A\sum_{k\in [1,n]\setminus K(\e,x)}1_B\circ T^kdm=o(a(n))\ \text{as}\ n\to\infty.$$
Now let $g\in C_B(\Bbb R_+)_+$, fix $\D>0$ and let $\e>0$ be so that
$$g(y)=g(m(f))\pm\D\ \forall\ y\in \Bbb R_+,\ |y-m(f)|< \e m(f).$$ It follows that
\begin{align*}\sum_{k=1}^n&|\int_Ag(\tfrac{s_k}{a(k)})1_B\circ T^kdm-g(m(f))m(A\cap T^{-k}B)|\\ &=
\sum_{k=1}^n|\int_{A\cap T^{-k}B}(g(\tfrac{s_k}{a(k)})-g(m(f))dm|\\ &
\le\int_A\(\sum_{k=1}^n1_B\circ T^k|g(\tfrac{s_k}{a(k)})-g(m(f))|\)dm.\end{align*}
Now
\begin{align*}\sum_{k=1}^n&1_B\circ T^k|g(\tfrac{s_k}{a(k)})-g(m(f))|\\ &=
\(\sum_{k\in  [1,n]\cap K(\e,x)}+\sum_{k\in  [1,n]\setminus K(\e,x)}\)1_B\circ T^k|g(\tfrac{s_k}{a(k)})-g(m(f))|\\ &=:I_n+II_n; 
\end{align*}
and
\begin{align*}
\int_A&II_ndm\\ &\le 2\|g\|_{C_B}\int_A\sum_{k\in [1,n]\setminus K(\e,x)}1_B\circ T^kdm\\ &=o(a(n))\ \text{as}\ n\to\infty; 
\end{align*}
whereas
$$I_n\le \D \sum_{k\in  [1,n]\cap K(\e,x)}1_B\circ T^k$$
and
$$\int_A I_ndm\le \D\sum_{k=1}^nm(A\cap T^{-k}B)\sim \D m(A)m(B)a(n).$$

\

\

\

\

\

This proves\ref{dschemical}.\ \Checkedbox

\

\demo{Proof of\ \  \ref{faTrain}}
\

Using {\tt RWM} $\&$\ref{dschemical},
\begin{align*}\sum_{k=1}^n&|\int_Ag(\tfrac{s_k}{a(k)})1_B\circ T^kdm-g(m(f))u(k)m(A)m(B)|\\ &\le 
\sum_{k=1}^n|\int_{A\cap T^{-k}B}(g(\tfrac{s_k}{a(k)})-g(m(f)))dm|+\\ &\ \ \ \ \ +|g(m(f))|
\sum_{k=1}^n|m(A\cap T^{-k}B)-u(k)m(A)m(B)|\\ & =o(a(n))
\ \text{as}\ n\to\infty.\ \CheckedBox 
\end{align*}

The following extends ({\scriptsize\faLightbulbO}) in \cite{A-Sera}
\proclaim{Corollary 7\ \ ({\scriptsize\rm discrete time tied-down renewal})}
\

\ \ Suppose that $(\xi_1,\xi_2,\dots)$ are $\Bbb N$-valued iidrvs so that 
$\ell(t):=\Bbb E(\xi\wedge t)$ is slowly varying, then with $s_n:=\sum_{k=1}^n\xi_k$, 
$$\tfrac1{a(N)}\sum_{n=1}^N|\sum_{k=1}^ng(\tfrac{k}{a(n)})P([s_k=n])-g(1)u(n)|\xrightarrow[n\to\infty]{}0\ \forall\  g\in C_B(\Bbb R_+)$$ where  $u(n):=\sum_{k=1}^nP([s_k=n])$ and $a(n):=\sum_{k=1}^nu(k)\sim \frac{n}{\ell(n)}$.
  \endproclaim\demo{Proof}\ \ Let $(X,m,T)$ be the corresponding renewal shift-- a pointwise dual ergodic measure preserving
transformation  with $a_n(T)\sim a(n)$.
  By the Darling Kac theorem \ref{faCab} (on p.\pageref{faCab}) holds and  the result follows from Proposition 6\ref{dschemical}.\ \ \CheckedBox
\section*{\S6  Appendix}

\subsection*{Mittag-Leffler processes and local time processes}
\

For $I=[0,R]$ or $[0,\infty)$, let

$$D(I):=\{\om:I\to\Bbb R:\ \forall\ t\in I\ \exists\ \om(t\pm)\ \&\ \om(t+)=\om(t)\}$$
where $\om(t\pm):=\lim_{s\to t\pm}\om(s)$ and $\om(0-):=\om(0)$:\   the space of all c\`adl\`ag paths from $I$ to $\Bbb R$ and write
$D_R:=D([0,R])\ \&\ D_\infty:=D([0,\infty))$.
\

For $0<\gamma\leq 1/2$, let $X\in\text{\tt RV}(D_\infty)$ be a symmetric $1/(1-\gamma)$-stable L\'evy process 
with $X_0=0$. 
\

As in   \cite[Chapter V]{Ber}, denote by $p_t(x)=\mathbb{P}[X_t\in dx]/dx$ its  transition density and set $\frak c:=2\int_0^1 p_s(0)ds$,  then $X$ admits the {\it local time process} at $0$, $L\in\text{\tt RV}(D_{\uparrow,\infty})$: 
\begin{align*}
	L(t)=\lim_{\varepsilon \to 0+}\frac{1}{\frak c\varepsilon}\int_0^t 1_{\textstyle \{|X_s|<\varepsilon\}}ds.
\end{align*}
By \cite[Proposition V.4]{Ber},  $L\ =\ \frak m_\g.$
\

As shown in \cite{Ber} (Theorem IV.4 and Proposition V.4), the closed support of the Stieltjes measure $\mu_L$ on $[0,\infty)$ defined by  $L$  is  $\overline{\{t\geq0:X(t)=0\}}$.
\

Thus a.s., the last visit time of $X$ at $0$ before time $1$ is given by
\begin{align*}G&:=\sup\{0\leq t\leq1: X_t=0\}=\sup([0,1]\cap {\rm supp}(\mu_L))\\ &
    =\mathcal{G}(L)(1)=\mathcal{G}(\frak m_\g)(1)
    \end{align*}
where $\mathcal{G}$ is as in $\S0$.
\

By \cite[Theorem VIII.12]{Ber}  the {\it symmetric $1/(1-\gamma)$-stable L\'evy bridge process}  (from $0$ to $0$ of length $1$) 
$X^{\rm( br)}\in\text{\tt RV}(D_1)$ is given by 
\begin{align*}
    X^{\rm( br)}(t)=\frac{X_{Gt}}{G^{1-\gamma}}.
\end{align*}

Moreover, $G$ is independent of $X^{\rm (br)}$ and is Beta$(\gamma, 1-\gamma)$-distributed.
\

As before, the local time process $L^{\rm(br)}\in \text{\tt RV}(D_{\uparrow,1})$ of $X^{\rm(br)}$ at $0$ satisfies
\begin{align*}
	L^{\rm( br)}(t)&\xleftarrow[\e\to 0+]{×}\frac{1}{\frak c\varepsilon}\int_0^t 1_{\textstyle \{|X^{\rm(br)}(s)|<\varepsilon\}}ds\\ &=
	\frac{1}{\frak c\varepsilon}\int_0^t 1_{\textstyle \{|X(Gs)|<G^{1-\g}\varepsilon\}}ds
	\\ &=
	\frac{1}{G\frak c\varepsilon}\int_0^{Gt} 1_{\textstyle \{|X(u)|<G^{1-\g}\varepsilon\}}du\ \ \text{changing variables}
	\\
	&\xrightarrow[\e\to 0+]{×}
	\frac{\frak m_\g(Gt)}{G^\gamma}=\frak w_\g(t),
\end{align*}
where $\frak w_\g$ the tied-down $\gamma$-ML process as in \S0.
\

The process $\frak w_\g=L^{\rm (br)}$ is independent of $G$, since $L^{\rm(br)}$ is a functional of $X^{\rm(br)}$.
\

Similarly, for $0<\gamma<1$,  the local time of the $(2-2\gamma)$-dimensional Bessel process is also  supported on the closure of the zero set and equal distributionally to $\frak m_\g$ (\cite{Mo-Os}). The $(2-2\gamma)$-dimensional Bessel bridge from $0$ to $0$ of length $1$ is constructed e.g.
in \cite[Theorem 3 and Section 2.2]{ChaUri}. The local time of this bridge process is shown to be $\frak w_\g$ in \cite{BPY}.\footnote{For more details about Bessel processes, bridges and their local time processes
see \cite{FPY, PitYor97, RevYor, RogWil} .}

\subsection*{A proof of Proposition C via Bessel bridges}
\

We'll deduce Proposition C as a consequence of \cite[Proposition 2]{FPY}. Following the notation there, for $r\in (0,\infty]$, let
$$D_{+,r}:=\{\om\in D_r:\ \om\ge 0\}.$$
Defining $X_t:D_{+,\infty}\to [0,\infty)$ by $X_t(\om):=\om(t)$, set $\mathcal{F}_{t-}=\sigma(X_s:0\leq s < t)$ and $\mathcal{F}_\infty=\sigma(X_s:s\geq0)$.  
\

Fix  $0<\gamma<1$ and for $x\ge 0$, let  $\mathbb{P}_x\in\mathcal P(D_{+,\infty})$ denote the law of a $(2-2\gamma)$-dimensional Bessel process on $[0,\infty)$ starting from $x$.
\

For $t>0$ and $x,y\ge 0$, 
the {\it $(x,t,y)$-bridge law} $\mathbb{P}_{x,y}^{t}$ on $(D_{+,\infty}, \mathcal{F}_{t-})$  is the regular conditional distribution of $\mathbb{P}_x|_{\mathcal{F_{t-}}}$ given $X_{t-}=y$ (see \cite[Section 2]{FPY} and \cite[Sections 1 and 2]{ChaUri}).

Denote the {\it usual augmentation} (as in \cite{RevYor}) by $\mathcal{F}^\ast=(\mathcal{F}_t^\ast:\ t\ge 0)$ which is the smallest right continuous filtration where each $\mathcal{F}_t^\ast$ is a completion of  $\sigma(X_s:0\leq s \leq t)$ with respect to $(\mathbb{P}_x: x\ge 0)$ and
let $(H_t:t\geq0)$ be a non-negative predictable process and $(A_t:t\geq0)$ be a continuous additive functional with respect to  $(\mathcal{F}_t^*:t\geq0)$. 
\

Assuming $\mathbb{P}_x[A_t<\infty]=1$, for any $x$ and $t$, \cite[Proposition 2]{FPY} yields that
\begin{align*}\tag*{\faFlask}\label{faFlask}
	\mathbb{E}_x\left[\int_0^t H_s dA_s\right]
	=
	\mathbb{E}_x\left[\int_0^t \mathbb{E}^{s}_{x, X_s}[H_s]dA_s\right],
\end{align*}
where $\mathbb{E}_x$ and $\mathbb{E}_{x,y}^s$ denote the expectations with respect to $\mathbb{P}_x$ and $\mathbb{P}_{x,y}^s$, respectively.
\

The local time at $0$: $L\in \text{\tt RV}(D_{\uparrow,\infty})$   with $\mathbb{E}_0[L_1]=1$ is a continuous additive functional and hence a predictable process. 
\

In addition we see $\int_0^t 1_{\{X_s=0\}}dL_s=L_t$ for any $t$. See \cite[Example]{FPY} or \cite[Chapter IV]{Ber}. 
\

Substituting $A_s=L_s$ into the \ref{faFlask}, we obtain
\begin{align*}\tag*{\faBeer}\label{faBeer}
   \mathbb{E}_x\left[\int_0^t H_s dL_s\right]
	=
	\mathbb{E}_x\left[\int_0^t \mathbb{E}^{s}_{x, 0}[H_s]dL_s\right]. 	
\end{align*}

Recall $\D_{s,\g}L=(L_{st}/s^\gamma:0\leq t\leq 1)$ and $L^{(1)}=(L_t:0\leq t \leq 1)$. Let $h:D_{\uparrow, 1}\to \mathbb{R}$ be a bounded measurable functional.
\

To prove Proposition C, apply  \ref{faBeer} with $x=0$, $t=1$ and $H_s=h(\D_{s,\g}L)$ and  obtain
                                                                                                          
\begin{align*}
\mathbb{E}_0\biggl[\int_0^1 & h(\D_{s,\g}L) dL_s\biggr]
	=
	\mathbb{E}_0\left[\int_0^1 \mathbb{E}_{0, 0}^s[h(\D_{s,\g}L)]dL_s\right]\ \text{by \ref{faBeer}}
	\\
	&=
	\mathbb{E}_0\left[\int_0^1 \mathbb{E}_{0, 0}^1[h(L^{(1)})]dL_s\right]\ 
	\because\ \mathbb{E}_{0, 0}^s[h(\D_{s,\g}L)]=\mathbb{E}_{0,0}^{1}[h(L^{(1)})]
	\\
	&=
	\mathbb{E}_0[L_1]\mathbb{E}_{0, 0}^1[h(L^{(1)})]=\mathbb{E}_{0, 0}^1[h(L^{(1)})]
	.
\end{align*}
This completes the proof, since  $L=\frak m_\g$ under $\mathbb{P}_0\ \&\ L^{(1)}=\frak w_\g$ under $\mathbb{P}_{0,0}^1$
 as explained above.\ \ \Checkedbox

\end{document}